\documentclass[11pt]{article}
\usepackage{authblk}
\usepackage[backend=bibtex,style=numeric-comp]{biblatex}

\title{Fitting ODE models of tear film breakup}
\author[1]{Tobin A. Driscoll}
\author[1]{Richard J. Braun}
\author[2]{Rayanne A. Luke}
\author[1]{Dominick Sinopoli}
\author[1]{Aashish Phatak}
\author[1]{Julianna Dorsch}
\author[3]{Carolyn G. Begley}
\author[4]{Deborah Awisi-Gyau}

\affil[1]{Department of Mathematical Sciences, University of Delaware, Newark, DE 19716, USA}
\affil[2]{Department of Applied Mathematics and Statistics, The Johns Hopkins University, Baltimore, MD 21218 USA}
\affil[3]{School of Optometry, Indiana University, Bloomington, IN 47405, USA}
\affil[4]{Alcon Research LLC, 6201 South Freeway, Fort Worth, TX 76134, USA}

\date{}

\usepackage[letterpaper,total={6.5in,8in}]{geometry}

\usepackage{amsmath,amssymb,amsthm}
\usepackage{siunitx}
\usepackage{appendix}
\usepackage{url}
\usepackage[]{graphicx}
\usepackage[hidelinks]{hyperref}
\graphicspath{ {.}{../figures} }
\usepackage{xcolor}
\hypersetup{
    colorlinks,
    linkcolor={red!50!black},
    citecolor={blue!50!black},
    urlcolor={blue!80!black}
}
\newcommand{\micron}{\ensuremath{\mu\text{m}}}
\newcommand{\mupermin}{\micro\meter/\minute}
\newcommand{\figref}[1]{\autoref{#1}}

\begin{document}

\maketitle

\begin{abstract}
 
\textit{Purpose.} Several elements are developed to quantitatively determine the contribution of different physical and chemical effects to tear breakup (TBU) in subjects with no self-reported history of dry eye or other ocular surface disease.  Fluorescence (FL) imaging is employed to visualize the tear film and to determine tear film (TF) thinning and potential TBU.

\textit{Methods.} An automated system using a convolutional neural network is deployed that was trained and tested on more than 50,000 images from FL imaging experiments.  The trained system could identify multiple TBU instances in each trial.  Once identified, extracted FL intensity data was fit by mathematical models that included tangential flow along the eye, evaporation, osmosis and FL intensity of emission from the tear film.  The mathematical models consisted of systems of ordinary differential equations for the aqueous layer thickness, osmolarity, and the FL concentration; they are a local approximation to TF thinning and/or TBU dynamics.  FL intensity was computed using the resulting thickness and FL concentration.    Optimizing the fit of the models to the FL intensity data determined the mechanism(s) driving each instance of TBU and produced an estimate of the osmolarity within TBU.

\textit{Results.} Initial estimates for FL concentration and initial TF thickness agree well with prior results. Fits were produced for $N=467$ instances of potential TBU from 15 non-DED subjects. The results showed a distribution of causes of TBU in these healthy subjects, as reflected by estimated flow and evaporation rates, which appear to agree well with previously published data. Final osmolarity depended strongly on the TBU mechanism, generally increasing with evaporation rate but complicated by the dependence on flow.

\textit{Conclusion.} The method has the potential to classify TBU instances based on the mechanism and dynamics and to estimate the final osmolarity at the TBU locus. The results suggest that it might be possible to classify individual subjects and provide a baseline for comparison and potential classification of dry eye disease subjects.

 \end{abstract}
 
\section{Introduction}
\label{sec:intro}

In this paper, we generate quantitative estimates of important parameters for the tear film on the surface of the eye in healthy subjects.  We do this with what we believe, at the time of writing, to be unprecedented precision and quantity.  The dataset creates a preliminary baseline for a small population of subjects without dry eye disease (DED).  The importance of this baseline is that it may be used to contrast what is found for a population with DED, thus leading to better understanding of the mechanisms at work in this disease that affects millions of people \cite{schaumbergPrevalanceWomen2003,schaumbergPrevalenceMen2009,stapletonBookepi2015,stapletonDEWSIIepi2017}. Though this work does not give a complete baseline for non-DED eyes, or a contrast with data for DED eyes, we develop the method in detail and explain how it can reveal the mechanisms behind individual instances of thinning and tear breakup (TBU) in the tear film (TF).

The introduction is structured as follows. Firstly, we give some background on the tear film, ocular surface and DED.  Secondly, we briefly discuss some related methods for imaging the tear film.  Thirdly, we discuss methods to extract data about tear film dynamics.  Finally, we discuss mathematical models for tear film dynamics, and best fits of those models to data extracted from the tear film.

\emph{Tear Film}  The TF plays an important role in vision and ocular surface health \cite{nelson2017tfos}.   The TF is established during a blink, and lubricates the cornea and the conjunctival surfaces lining the gap between the lids and the globe \cite{pultSpontaneousBlink2015}. The air/tear film interface causes the tear film to have the most powerful refractive surface in the eye; thus, keeping that surface smooth and regular is essential to clear vision \cite{ThibosClinicalApplications1999}. When the TF fails to uniformly coat the ocular surface, it is said that tear breakup has occurred \cite{nornDesiccationPrecorneal1969,choTearBreakup1998}. TBU may cause the ocular surface to be exposed to cooling \cite{Mapstone68a,EfronYoung89,Belmonte_Cold_2011} and evaporation \cite{nicholsThinningRatePrecorneal2005,DurschFLandThermal2017}, and evaporation may lead to tear hyperosmolarity \cite{DEWSdef,liuLinkInstability2009,BraunKing-Smith2015,craigTFOSDEWSdef2017} and mechanical stimulus to the surface \cite{awisi-gyauChangesCorneal2019}. The exposure of the ocular surface to hyperosmolarity from TBU is thought to play a central role in the etiology of DED \cite{DEWSdef,craigTFOSDEWSdef2017} which affects millions of people \cite{stapletonDEWSIIepi2017}.  As a result of this significance, TBU dynamics have been studied for more than 50 years using a variety of methods \cite{nornDesiccationPrecorneal1969,wolffsohnTFOSDEWS2017}.  Clinically, the instability of the tear film is measured by the technique of tear breakup time (TBUT), in which the time to the first break or irregularity of the tear film is measured.

\emph{Imaging methods}  The imaging methods for TBU dynamics are numerous.  Here we list a few of them: visualization with dyes such as fluorescein (FL) \cite{nornDesiccationPrecorneal1969,choReliabilityTear1992}; reflection of a pattern using a grid \cite{mengherProposedNIBUT1985} or placido disc images \cite{llorensVideokeratoscopy2019}; interferometry and spectrometry \cite{Doane1989,danjoInterferometry1994,gotoKineticAnalysis2003,King-SmithFink2004,SegevDynamicAssessment2020}; simultaneous imaging with fluorescence (FL) imaging and retroillumination  \cite{BraunKing-Smith2015};  and simultaneous FL imaging with interferometry \cite{King-SmithIOVS13a}.

These and other approaches have quantified various aspects of TF parameters such as thicknesses, thinning rates, TBUTs and more.  In this work, we focus on fluorescence imaging as an experimental method to collect data on aqueous layer (AL) dynamics. This method is chosen due to the relatively low cost, ease of use and widespread use in the clinic.  Clinically, short TBUTs indicate an unstable TF and the possible presence of DED \cite{wolffsohnTFOSDEWS2017}.
Despite the utility of the method, repeatability from one clinician or researcher to the next and one clinic to the next can be a challenge \cite{nornTearFilm1986}, though some maintain that TBUT measurements can be generally repeatable under some circumstances \cite{choReliabilityTear1992}.  In this work, we aim to use automated detection of FL imaging to (i) repeatably extract FL imaging data of TF thinning and TBU, and subsequently to (ii) optimize the fit of mathematical models to that data to identify mechanism and (iii) estimate important parameters within TBU.

Efforts to automate TBU and DED measurements were recently reviewed by Vyas and Mehta \cite{vyasComprehensiveSurvey2020}. Early efforts generally aimed at quantifying TF breakup time measurement and related quantities\cite{yedidyaAutoTBUSequenceDetect2008,ramosAnalysisParameters2014}.
Vyas and Mehta \cite{vyasComprehensiveSurvey2020} surveyed various methods for automating measurements and diagnoses, including: tear meniscus evaluation using optical coherence tomography \cite{bartuzelTearMeniscusMeas2014};  thermal imaging to attempt to diagnose DED \cite{acharyaAutomatedDEDDiagThermal2015}; and fluorescence imaging of the TF for tear breakup time detection \cite{suTearFilm2018} and DED diagnosis \cite{remeseiroCASDESComputerAided2016}.

\emph{Extraction of data}  Our method in this paper is adapted from that of Su et al \cite{suTearFilm2018}. In their system, a convolutional neural network (CNN) is implemented that determines a region of interest where TBU is most likely to occur.  Then, the region of interest is followed in time and the first frame where TBU is found determines the TBUT.  Their method is trained on TBU and TBUT data from experienced clinical researchers, and is therefore designed to imitate the clinical determination of TBUT for the purpose of DED diagnosis.
While we retained the CNN design from their work, we introduced several changes to the approach of Su et al \cite{suTearFilm2018}. The method is adapted to identify multiple regions of TBU in every trial.  We extracted a time series of FL thinning data from each TBU region.  We used that FL imaging time series to determine TBUT (if appropriate) as well as optimal parameters for mathematical models to determine important quantities of interest with thinning and TBU areas.  The optimal parameters allow us to identify the mechanism(s) driving each instance of TBU.

\emph{Mathematical models}  A variety of mathematical modeling approaches for the TF have been developed.
For overall flows and concentrations of interest in the TF, there have been compartment models, systems of ordinary differential equations (ODEs), or differential algebraic equations (DAEs) that have included the effect of blinks
\cite{gaffneyTFOsmolarityModel2010,cerretaniTFDynamicsModel2014} and contact lenses
\cite{gauseMechModelingDrugDel2016,kimContactLensesForProtect2022}.   TBU and TF dynamics with contact lenses are beyond the scope of this paper.

A few categories of 1D partial differential equation (PDE) models in space and time have been developed; this includes
TF drainage for the open eye during the interblink \cite{SharmaTiwari98,WongFatt1996,MillerPolse02}. Those models used a Newtonian fluid close to water in viscosity and measured TF values.  Boundary conditions (BCs) at the end of the film mimicked the TF and drove flows to redistribute TF.  Effects added to this type of model include Marangoni effects \cite{BergerCorrsin74}, evaporation \cite{BraunFitt03}, van der Waals wetting terms \cite{WinterAnderson10} and curvature of the ocular surface \cite{braunEllipsoidalSubstrate2012}.  Local models for TF thinning and TBU include those which have been studied for the following effects:  evaporation to air and osmosis from corneal surface \cite{PengEtal2014} and with
fluorescence \cite{braunTearFilm2018};   Marangoni effects \cite{zhongFLimaging2019};  a non-polar lipid layer (LL) \cite{brunaInfluenceNonpolar2014,stapfDuplexTearFilm2017}; dewetting of the ocular surface from long-range van der Waals forces \cite{sharmaMechTFRupture1985,sharmaLongRange1999}; dewetting of the ocular surface with mucin-dependent viscosity \cite{deyContinuousMucinProfiles2019,deyContinuousMucinCorrection2020} and membrane-associated mucins \cite{choudhuryMembraneMucin2021}.  Some models for TBU are discussed in more detail below.

Models for TF formation, which occurs during the opening phase of the blink cycle, have been studied as well.  A seminal work in this area is Wong et al \cite{WongFatt1996}, which treated the TF deposition as a thin film coating flow model; this is a cornerstone of later papers although they modified the approach. Later models have included the effect of polar lipids via the Marangoni effect \cite{JonesPlease2005,JonesMcElwain2006,aydemirEffectPolar2011,makiLipidReservoir2020}; partial blinks \cite{heryudonoSingleequationModels2007,dengPartialBlinks2013};  a non-polar LL \cite{brunaInfluenceNonpolar2014,ZubkovBreward12}; the curvature of the ocular surface \cite{allouche2017} and non-Newtonian effects \cite{JossicLefevre2009,mehdaouiShearThinning2021,mehdaouiGelBased2021}.

Models for flow over the (2D) exposed ocular surface have been developed \cite{makiTearFilm2010a,liOsmoDynamics2016,LiBraun17,broschSimulationThin2016,braunMathematicalModels2019,makiModelTear2019}.  The 2D models capture a number of aspects of the overall flows, osmolarity and fluorescence imaging. Some 2D models may take into account the effect of blinking via time-dependent flow BCs with no lid motion \cite{LiBraun14}, or via lid motion with model problems plus simple BCs \cite{broschSimulationThin2016}, but there is much room to develop blinking models.

Local models have been developed for flow in TBU regions.  
Peng et al.\cite{PengEtal2014} studied TBU driven by tear evaporation through a LL distribution that was fixed in space. In their model, evaporation rate depended on the temperature of the ocular surface, as well as the temperature, relative humidity and wind conditions of the surroundings.  They found that evaporation could drive the AL thickness to very small values and  thus TBU.  Simple ODE models of TF thinning with osmosis could develop sufficiently elevated osmolarity that could stop thinning and TBU \cite{braunDynamicsTear2012,BraunKing-Smith2015}; however, Peng et al.\cite{PengEtal2014} found that diffusion of osmolarity (salt ions) out of the high concentration region within TBU prevented sufficient osmosis to stop TBU \cite{PengEtal2014}.   A dynamic LL was introduced in Stapf et al.\cite{stapfDuplexTearFilm2017}.  The model consisted of two Newtonian layers: a relatively thick and less viscous shear layer topped by a relatively thin but more viscous extensional layer through which evaporation occurred.  Stapf et al.\cite{stapfDuplexTearFilm2017} found that TBU could occur, but the model could yield longer TBUTs than would be observed \emph{in vivo}.  This also happened with models that incorporated mucin effects \cite{deyContinuousMucinProfiles2019,choudhuryMembraneMucin2021}.

Braun et al. \cite{braunTearFilm2018} simplified TF dynamics to a single layer for the AL with evaporation modeled as a fixed Gaussian, but they included fluorescein concentration and fluorescence in their models of TBU.  They found that the fluorescence dynamics depended on initial FL concentration, evaporation distribution width (related to TBU size) and film thickness in a complicated way, but the mechanisms at work in various instances were clarified by the model.  Subsequently, models were proposed to include rapid thinning that could be induced by excess lipid acting as a surfactant \cite{zhongMathematicalModelling2018,zhongFLimaging2019,lukeParameterEstimationMixedMechanism2021a}.  The models explained many aspects of TBU, but they tend to overestimate the size of the TBU region \cite{lukeParameterEstimationMixedMechanism2021a}.

In this work, we use local models for tear break up involving tangential flow, evaporation, osmosis and fluorescence, but the models have been simplified to ODEs for the thickness, osmolarity, fluorescein and fluorescent intensity\cite{lukeParameterEstimationMixedMechanism2021b}.  We find the optimal parameters for these models that make them as close as possible to FL intensity data extracted from video recordings of \emph{in vivo} TFs.  With those optimal parameters, we can infer which effects were most important in each TBU instance.  We use a CNN to extract data for many TBU instances in order to get a more complete picture of TBU for the cohort of healthy subjects studied.

\emph{Paper structure}  This paper is structured as follows.  The methods section will describe in some detail the FL imaging used to generate data; the extraction method we used to obtain the detailed thinning data; and mathematical methods and models used to fit that data and determine TBU parameters of interest.  In the results section, we present the results of applying these methods.  In the discussion section we explain the context and significance of the results. In the conclusion section, we summarize our findings and discuss possible future directions.

\section{Methods}

\subsection{Fluorescence imaging}

The experimental data was collected at Indiana University and was approved by the Biomedical Institutional Review Board of Indiana University. The principles of the Declaration of Helsinki were followed during data collection, and informed consent was obtained from all subjects. Data collection is described in a previous publication~\cite{awisi-gyauChangesCorneal2019} and discussed in several papers~\cite{awisi-gyauChangesCorneal2019,zhongFLimaging2019,lukeParameterEstimationMixedMechanism2021a,lukeParameterEstimationMixedMechanism2021b,lukeParameterEstimation2020}, but will be summarized briefly here.  Twenty-five subjects with no self-reported history of DED, ocular surface or systemic disease, ocular surgery or medications affecting ocular sensation participated in the study. Subjects were seated behind a slit lamp biomicroscope and 2 $\mu$l of 2\% sodium fluorescein solution was instilled in the subject’s eye. Subjects were asked to keep the tested eye open as long as possible (STARE trial) while the tear film was imaged with a cobalt blue excitation filter over the illumination system and a Wratten \#12 filter over the observation port.  With this illumination system, the aqueous layer of the TF fluoresced green~\cite{carlsonClinicalProcedures2004} with dark areas appearing due to TBU.

A trial is the sequence of images of the subject’s eye following a
few quick blinks. The trial records the fluorescence of the aqueous part of the TF. The trials typically start with an FL concentration close to 0.2\% (discussed more below), which is the so-called critical concentration where peak fluorescence occurs for thin TFs \cite{webberFLImaging1986}. The critical FL concentration may also be expressed as 0.0053 M \cite{lukeParameterEstimationMixedMechanism2021b}.

\subsection{TBU Detection}
\label{sec:tbu_detect}

We implemented a deep CNN \cite{HeDeepResidual2016,highamDeepLearningIntroduction2019} similar to the one used by Su et al.\cite{suTearFilm2018} to classify small square patches within an image as belonging to eyelids, eyelashes, sclera, TBU, and non-TBU. The architecture of the CNN is described in Table~\ref{tab:cnn} and requires a total of 313,637 parameters for training.

\begin{table}
    \centering
    \caption{Architecture of the neural network trained to classify $96\times 96$ RGB image tiles.}
    \begin{tabular}{lcccccl}
        Layer type & Number & Size & Stride, Pad & Output size & Activation \\
        Convolution & 32 & $5\times5$ & 1,2 & $96\times96\times32$ & ReLU \\
        Max pool & & $2\times2$ & 2,0 & $48\times48\times32$\\
        Convolution & 32 & $5\times5$ & 1,2 & $48\times48\times32$ & ReLU \\
        Average pool & & $3\times3$ & 2,1 & $24\times24\times32$\\
        Convolution & 64 & $5\times5$ & 1,2 & $24\times24\times64$ & ReLU  \\
        Average pool & & $3\times3$ & 2,1 & $12\times 12 \times64$ \\
        Convolution & 64 & $5\times5$ & 1,0 & $8 \times 8 \times64$ & ReLU \\
        Average pool & & $3\times3$ & 2,1 & $4 \times 4 \times64$ \\
        Convolution & 64 & $4\times 4$ & 1,0 & $1 \times 1 \times 128$ & ReLU \\
         Dropout, $p=0.4$ & & & & & & \\
         Dense & & & & 5 & softmax
    \end{tabular}
    \label{tab:cnn}
\end{table}

In order to obtain training data, we selected videos from 56 different trials spanning 10 distinct subjects, dividing still frames into tiles of $192 \times 192$ pixels (at 5.8 $\mu$m per pixel, tiles are 1.11 mm on a side). Selected tiles were then manually labeled into categories eyelid, eyelash, sclera, and TF, according to their most dominant feature in the context of the full image. The TF tiles were then sorted according to whether sufficiently many pixels were at luminous intensity 60 or less (on a scale of 0--255). The dark TF tiles were labeled as TBU, while the others were labeled as non-TBU. The total number of tiles labeled manually within each category is shown in Table~\ref{tab:tiles}. In order to accommodate the use of $96 \times 96$ tiles by the CNN, the labeled tiles were downsampled by a factor of two in each dimension. The data were then split so that 80\% (41838) were used for training and 20\% (10460) for testing.  The training set was artificially augmented by applying random flips and 90-degree rotations to the original tiles.

\begin{table}
    \centering
    \caption{Number of labeled tiles within each category for training the neural network.}
    \begin{tabular}{lccc}
        Category & Train & Test & Total \\ \hline
        Eyelash & 6460  & 1615 & 8075 \\
        Eyelid  & 2606  & 652  & 3258 \\
        non-TBU & 12976 & 3244 & 16220 \\
        Sclera  & 7879  & 1970 & 9849 \\
        TBU     & 11917 & 2979 & 14896 \\
        \hline
        Total & 41838 & 10460 & 52298
    \end{tabular}
    \label{tab:tiles}
\end{table}

The classification results from the test images are given in the confusion matrix shown in Table~\ref{tab:confusion}.  Of particular interest is the precision and recall for the TBU (0.93 and 0.96, respectively) and for the non-TBU (0.96 for both).  These results are more than adequate for our purposes of extracting thinning and TBU data.

\begin{table}
    \centering
    \caption{Confusion matrix showing the results from training the neural network.}
    \begin{tabular}{l|ccccc|r}
        Category & Eyelash  & Eyelid & non-TBU & Sclera & TBU  & Total \\ \hline
        Eyelash & 1400 & 98 & 17 & 32 & 68 & 1615\\
        Eyelid  & 112 & 494 & 2 & 7 & 37   & 652  \\
        non-TBU & 5 & 0 &  3108 & 70 & 61  & 3244  \\
        Sclera  & 38 & 20 & 72 & 1790 & 50 & 1970  \\
        TBU     & 22 & 15 & 42 & 53 & 2847 & 2979 \\
        \hline
        Total   & 1577 & 627 & 3241 & 1952 & 3063 & 10460 
    \end{tabular}
    \label{tab:confusion}
\end{table}

In order to detect TBU regions of interest (ROI), each trial video was first stabilized using the location of the Purkinje image of the lamp. Every frame of the stabilized video was overlaid with a grid of 192x192 overlapping pixel tiles with a stride length of 32 pixels. The tiles intersecting a detected corneal circle~\cite{driscollAutomaticDetectionCornea2021} were downsampled, and tiles marked by the CNN as very likely to be TBU were clustered to become ROIs. The locations were recorded relative to a rectangle cropped closely to the detected corneal circle. This process was continued throughout the video until at least three and as many as five distinct ROIs were identified. Additional details are given in Appendix C.

\subsection{Time series extraction}
\label{sec:data_extract}

Within each ROI, the images at each time were downsampled and subjected to a slight Gaussian blur, and a location was chosen within the ROI to sample the pixel intensity of the blurred image; details appear in Appendix C.

\begin{figure}
    \centering
    \includegraphics[width=\textwidth]{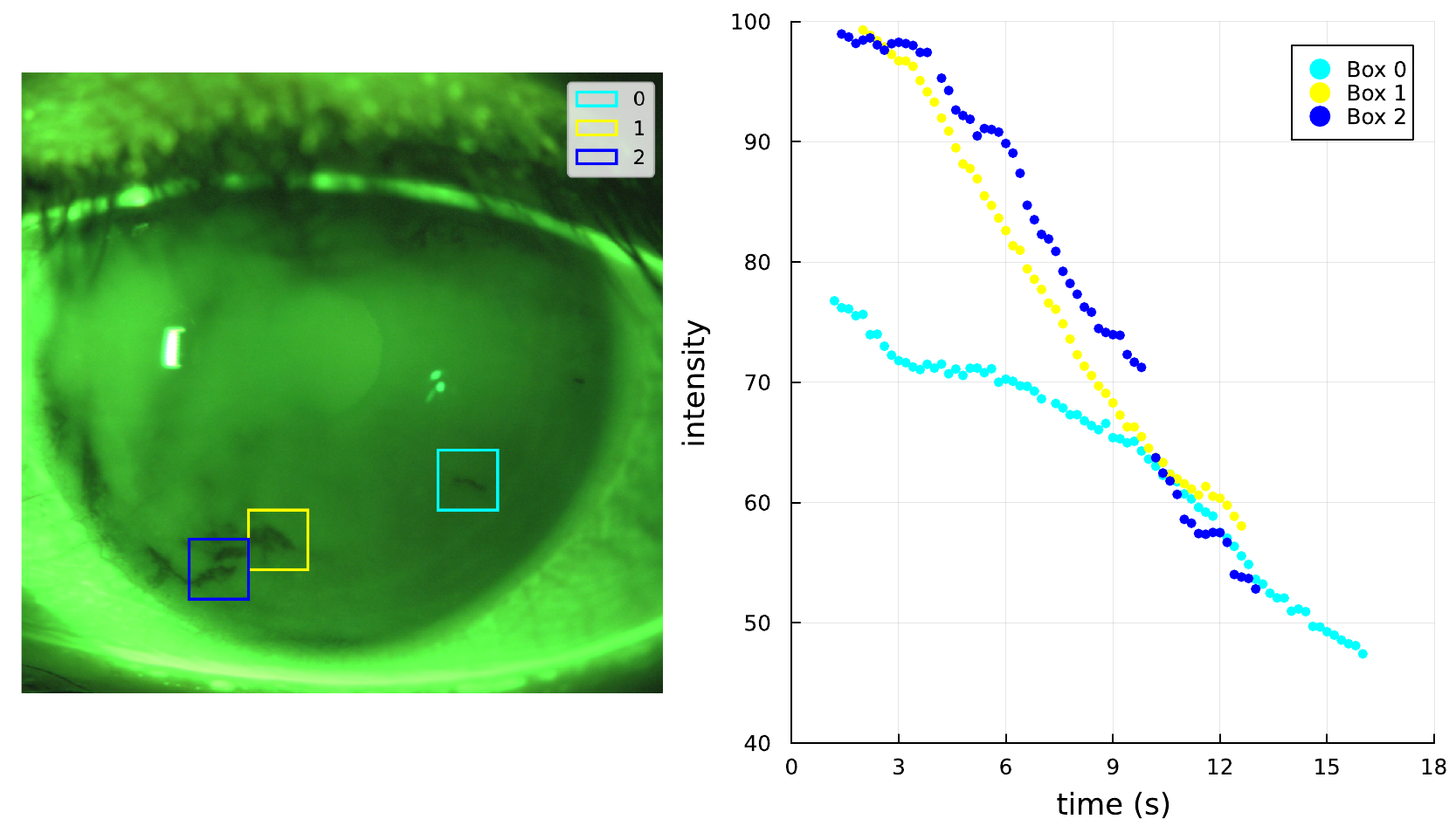}
    \caption{At left is one frame from the video of a trial, showing fluorescent intensity in green and the locations of likely TBU instances as boxes. Right shows the intensity time series captured for the marked boxes.}
    \label{fig:show_boxes}
\end{figure}

\figref{fig:show_boxes} shows example results of FL intensity data extraction for a single trial. The top left shows an image from late in the trial with likely TBU boxes marked, while the other plot shows intensity time series from the identified ROIs. As can be seen from the plots, the shape of the intensity curve can vary from one TBU instance to the next, even within the same trial. This phenomenon is not unexpected, since examples of different TBU mechanisms from the same subject have been reported in, e.g., simultaneous imaging experiments \cite{King-SmithIOVS13a}.
\clearpage

\subsection{Model fitting}
\label{sec:model-fitting}

\subsubsection{Models}
\label{sec:modelsIntro}

\begin{figure}[htb]
    \centering
    \includegraphics[width=4in]{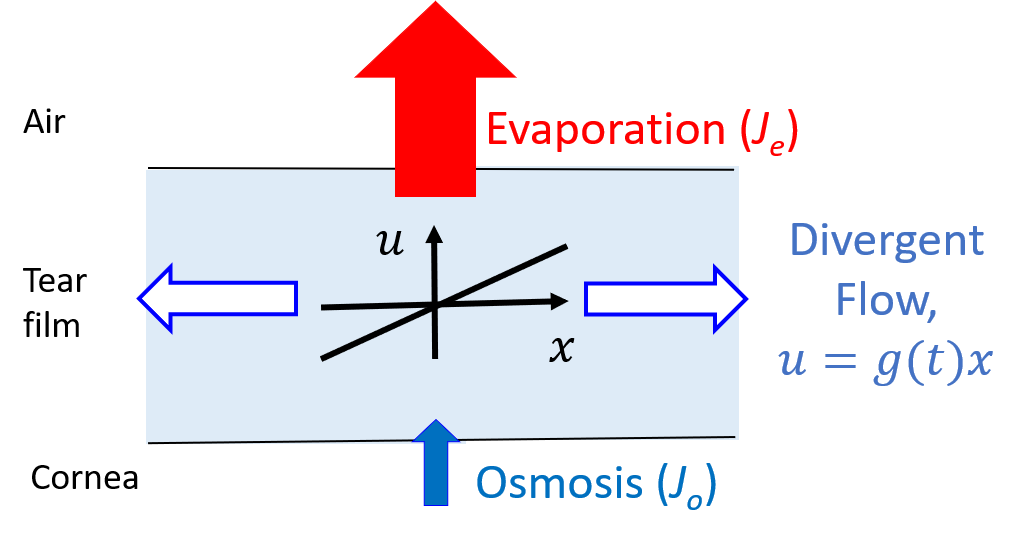}
    \caption{Sketch of ingredients in ODE models.  The film is spatially uniform.  It may be subject to loss of water via evaporation, supply of water due to osmosis, and to divergent flow away from the middle of the film. The thickness is given by $h(t)$.
    }
    \label{fig:model_sketch}
\end{figure}

A sketch showing the ingredients of the non-dimensional model are shown in \figref{fig:model_sketch}.
Evaporative loss of water is given, in dimensionless form, by the constant $J_e = v$.  The dimensionless supply of water from osmosis, which may result from hypertonicity due to evaporation, is given by $J_o = P_c(c-1)$. Here $P_c$ is a dimensionless permeability and $c$ is the
osmolarity.  The divergent flow is given by the velocity field $u=g(t)x$, and the strain rate $\partial u/\partial x = g(t)$ characterizes the flow.  This flow is constant throughout out the thickness of the film, but varies along its length; the fluid is simply being stretched.

We model TBU using a hierarchy of ODE models\cite{braunMathematicalModels2019,lukeParameterEstimationMixedMechanism2021b} represented as a system of nondimensional equations:
\begin{align}
    \label{eq:ode-h}
    \frac{dh}{dt} &= -g(t)h + P_c(c - 1) - v, \\
    \label{eq:ode-t}
    \frac{d(hc)}{dt} & = -g(t)hc,
\end{align}
for $0\le t \le 1$, after rescaling time by a scale $t_s$. Here the unknowns are $h(t)$, the TF thickness, and $c(t)$, the osmolarity.
The dependent variables are normalized so that $h(0)=c(0)=1.$ These values are found by scaling the dimensional variables (primed) with
\begin{equation}
    \label{eq:hc_scales}
    h = \frac{h'}{h_0}, \quad c = \frac{c'}{c_0},
\end{equation}
where $h_0$ is the initial film thickness and $c_0$ is the isotonic osmolarity.

The function $g(t)$ accounts for transverse flow and may take one of the following three functional forms:
\begin{subequations}
\label{eq:strain}
\begin{align}
    \textbf{Model O:} & \quad g(t) \equiv 0,\label{eq:strain-O} \\
    \textbf{Model F:} & \quad g(t) = a,  \label{eq:strain-F}  \\
    \textbf{Model D:} & \quad g(t) = b_1e^{-b_2t}. \label{eq:strain-D}
\end{align}
\end{subequations}
The values $a$, $b_1$, and $b_2$ are considered constant parameters. In Model O there is no fluid flow, so the model incorporates only evaporation and osmolarity. Model F adds constant extensional flow, while Model D allows extensional flow that decays to zero. Note that each model in \eqref{eq:strain} is a generalization of the models above it. We also considered an additional generalization allowing extensional flow that decays from one nonzero value to another ($g(t)=a+b_1e^{-b_2 t}$), but we do not report corresponding results due to relatively poor identifiability of its parameters for some of the data.

The parameters $v$, $a$, $b_1$, and $b_2$ (to the extent present) completely specify a model, and the dimensional versions are optimized to fit experimental FL intensity data. The parameters are nondimensional and related to their dimensional (primed) counterparts by
\begin{equation}
    \label{eq:dimensional}
    v = \frac{t_s v'}{h_0}, \quad  a = t_s a', \quad  b_1 = t_s b_1', \quad  b_2 = t_s b_2',
\end{equation}
where $t_s$ and $h_0$ are characteristic time and length scales, respectively. In practice, we choose $t_s$ as the duration of the observation window and $h_0$ as the initial thickness of the TF.  We also have $P_c = (P_o V_w c_0)/(h_0/t_s)$, where the dimensional permeability of the corneal surface, $P_o$, is fixed, $V_w$ is the molar volume of water, and $c_0$ is the isotonic osmolarity.  Values for these dimensional parameters are given in the appendix.  The permeability $P_o$ is not a parameter in the optimization because it is fixed\cite{BraunKing-Smith2015}; however, $P_c$ can vary while the other parameters are optimized.  Parameter values are given in Appendix B.

This spatially-uniform model allows us to obtain the nondimensional FL concentration $f(t)$ from scaling the osmolarity via
\begin{gather}
    \label{eq:FL_conc}
    f(t) = f_0 c(t),
\end{gather}
and the FL intensity from the film thickness and FL concentration via
\begin{gather}
    \label{eq:intensity}
    I(t) = I_0\frac{1 - \exp [ -\phi h(t)f(t) ]}{1 + f(t)^2},
\end{gather}
where $I$ is FL intensity and $\phi$ is the (nondimensional) Napierian extinction coefficient \cite{nichols2012,braunFluorescence2014}. Similarly to $P_c$, we have
$ \phi = \epsilon_f h_0 f_{cr}$, which includes the dimensional extinction coefficient $\epsilon_f$ (value given in Appendix B).  The value of $\phi$ varies from trial to trial because $h_0$ does.  The constant $I_0$ is used to match the initial observed intensity in the experiment.  Scaling both the experimental and theoretical FL intensities to start with unit value is desirable for the fitting to be described next.

\subsubsection{Fitting}
\label{eq:fitting}

For each video recording, the procedure of Wu et al.\cite{wuIncreasingOSStimulation2015} was used to estimate initial film thickness $h_0$ and initial fluorescein concentration $f_0$. We excluded as unreasonable all cases for which $h_0$ is outside the range \qtyrange{1}{10}{\micro\meter} or $f_0>0.35\%$. We use 
$f_0$ as the ratio between nondimensional osmolarity and FL concentration throughout the fit. The permeability parameter in~\eqref{eq:ode-h}, $P_c$, varies during the optimization as discussed in the previous section (see also Luke et al. (2021) \cite{lukeParameterEstimationMixedMechanism2021b}).

We excluded any intensity time series that showed substantial, sustained brightening; while this may happen \emph{in vivo} \cite{King-SmithIOVS13a}, we aim to fit thinning and TBU processes. Within each time series, instantaneous values that were local outliers were removed, and the time series was smoothed using an averaging filter.  An iterative procedure was then used to isolate a window of steepest average decrease lasting at least 3 seconds and excluding initial increases and final increases or plateaus. This window was judged to find the regime of thinning that the ODE models are best able to explain.
The intensity values were normalized by the initial value so that $I(0)=1$ and we determined $I_0$ in~\eqref{eq:intensity} to do that.

Given the normalized time series $I_k$ at times $t_0=0,t_1,\ldots,t_N=1$, the objective function for fitting ODE parameters was defined as the sum of squares,
\begin{align}
    \frac{1}{N} \sum_{k=0}^N \left[ I(t_k) - I_k \right]^2 = \frac{1}{N} ||I(t_k)-I_k||_2^2,
\end{align}
where $I(t)$ is from~\eqref{eq:intensity}, using the solution of~\eqref{eq:ode-h}--\eqref{eq:ode-t} and \eqref{eq:FL_conc} for $h(t)$ and $f(t)$. This objective was minimized over evaporation rate $v$ and the constants $a,\ b_1,\ b_2$ available in whatever form is chosen for $g(t)$. Constrained minimization was performed using both the BFGS and Nelder--Mead algorithms to confirm that the same minima were reached. The dimensional forms of the parameters were constrained to physically plausible ranges: $v'$ from \qtyrange{0}{40}{\micro\meter/\minute}, $a'$ from \qtyrange{-1}{2}{\per\second}, $b_1'$ from \qtyrange{-1}{5}{\per\second}, and $b_2$ from \qtyrange{0}{2}{\per\second}.

If, during ODE solution at a particular set of parameters, the numerical solution satisfied the conditions $\dot{h}(t)>0$ or $\dot{I}(t)>0$ over a sustained time interval, the solver was interrupted, and the optimization was given a penalty value to force selection of different values.  We made this choice because the models were designed for thinning of the TF.  Each optimization was attempted from multiple initializations in order to explore the global parameter space.

\begin{figure}
    \centering
    \includegraphics[width=\textwidth]{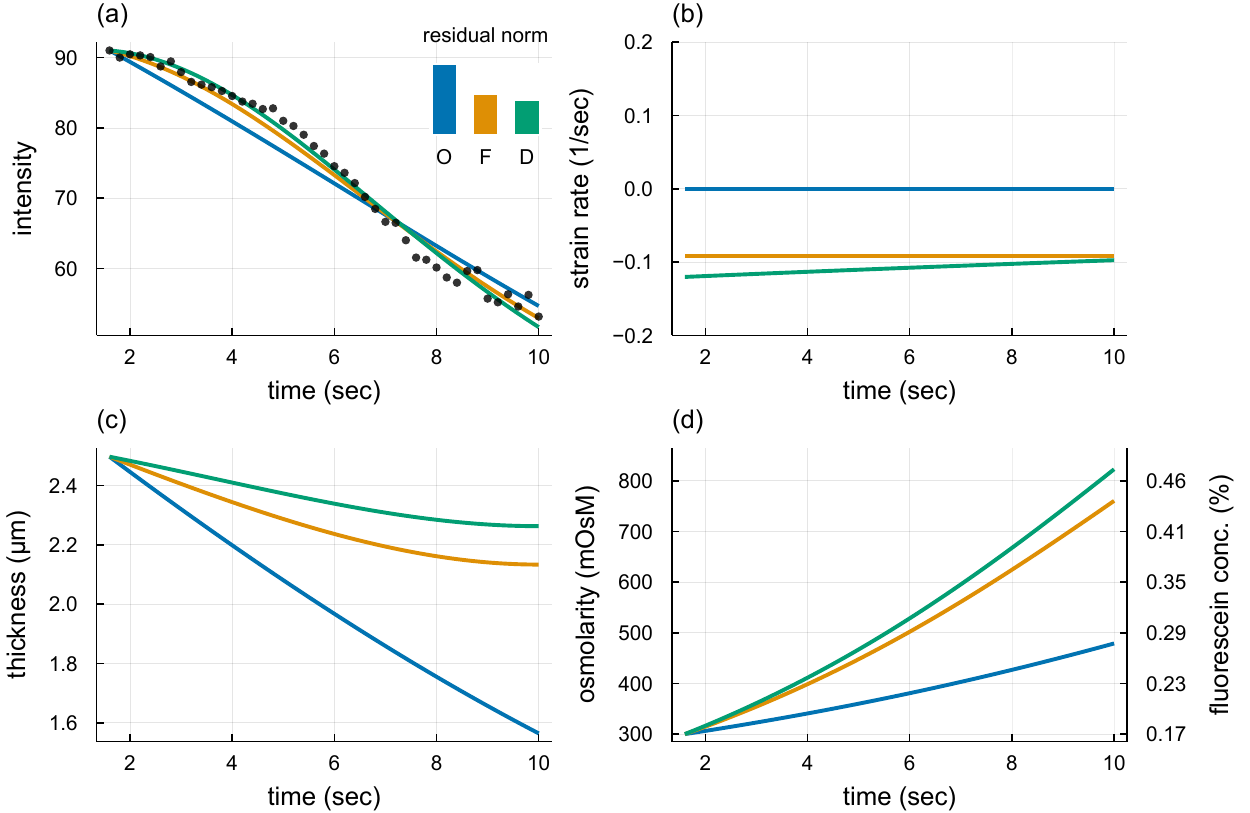}
    \caption{Results of fitting to a ROI with evaporation-dominated thinning. (a) FL intensity time series data (dots) and the best fits of the model types O, F, and D. The bar graph shows the relative residual norms of the fits. (b) The strain rate $g(t)$, showing a convergent flow in the models that allow it. (c) TF thickness in the three models. (d) Osmolarity and fluorescein concentration.
    }
    \label{fig:fits_evap_dominated}
\end{figure}

\begin{figure}
    \centering
    \includegraphics[width=\textwidth]{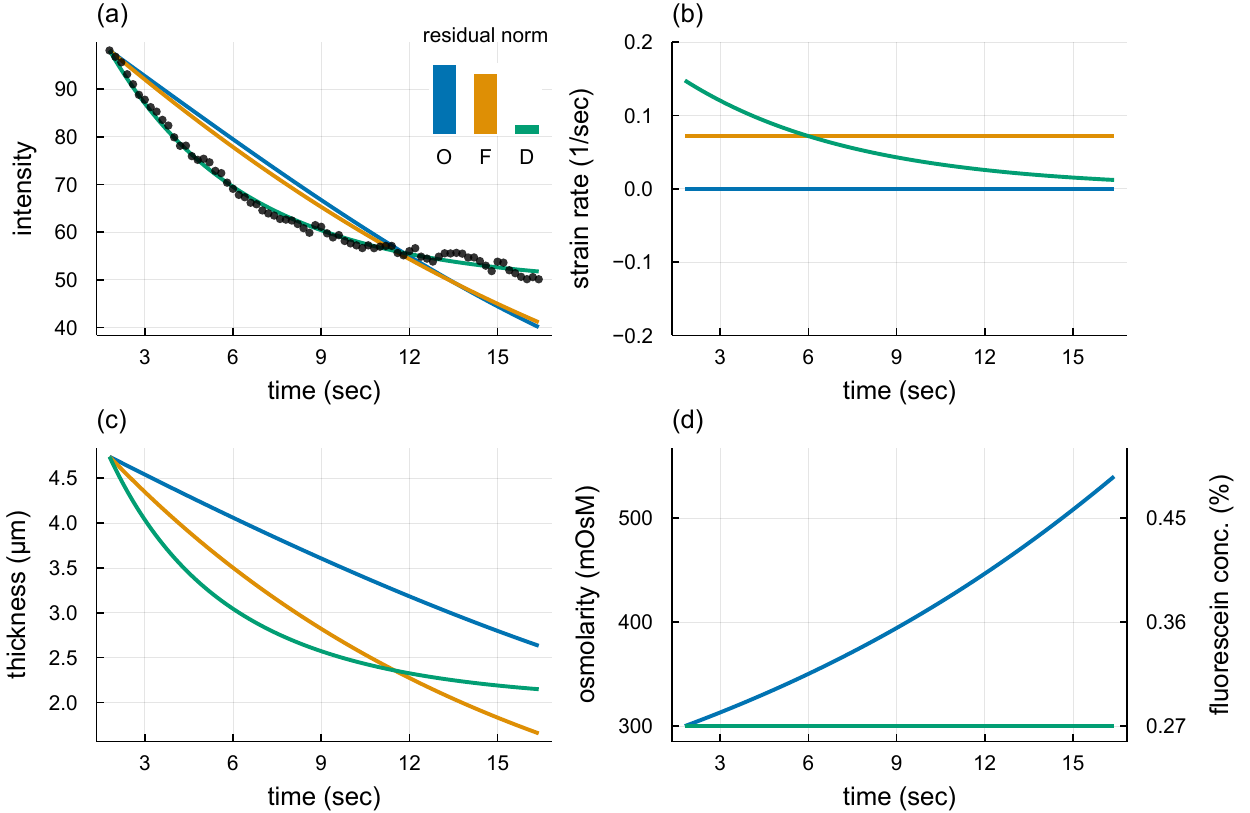}
    \caption{Results of fitting to a ROI with flow-dominated thinning. (a) FL intensity time series data (dots) and the best fits of the model types O, F, and D. The bar graph shows the relative residual norms of the fits. (b) The strain rate $g(t)$, showing a divergent flow in the models that allow it. (c) TF thickness. (d) Osmolarity and fluorescein concentration.
    }
    \label{fig:fits_flow_dominated}
\end{figure}

\figref{fig:fits_evap_dominated} shows fitting results for a particular ROI. \figref{fig:fits_evap_dominated}(a) shows that models F and D fit the data much better than does the evaporation-only model O. \figref{fig:fits_evap_dominated}(b) shows that the better models incorporate convergent flow to replace fluid lost to evaporation, which moderates the thinning (\figref{fig:fits_evap_dominated}(c)) but increases the osmolarity (\figref{fig:fits_evap_dominated}(d)). While model O found an optimal $v'$ at \qty{7.72}{\mupermin}, models F and D found $v'=17.8$ and \qty{20.0}{\mupermin}, respectively, indicating the dominance of evaporation in the thinning.

\figref{fig:fits_flow_dominated} shows fitting results for a different ROI. Here, model D is clearly superior, allowing a significant initial divergent flow that decays away. In this case, model O found $v'=$\qty{10.0}{\mupermin}, while the other models found $v'\approx 0$.  The osmolarity barely increases at all when flow is active (F or D), in contrast to the evaporative case.  The fluorescein concentration barely budges as well (proportional to the osmolarity), so the intensity change is due almost exclusively to the change in thickness (see (\ref{eq:intensity})).

Because the models form a hierarchy, the final residuals of the models must satisfy Model O $\ge$ Model F $\ge$ Model D.  As a result, the optimization of Model D was always best, and so results below are reported in terms of its parameters.  Those parameter values can change dramatically between different instances of thinning and/or TBU.

\section{Results}
\label{sec:results}

In total, 467 time series were successfully fitted to mathematical models. We begin discussing those cases by examining the initial conditions found for the analysis.

\subsection{Initial conditions for fitting}

We compare the distribution from the current results with other mathematical models and direct measurements of TF thickness.   Creech et al. \cite{creechTearFilmThickness1998} used the coating flow model of Wong et al. \cite{WongFatt1996} to estimate the thickness of the deposited TF from the opening phase of the blink.
Our initial estimates for thickness are close to those of published measurements \cite{nicholsThinningRatePrecorneal2005,wangTFThicknessOCT2003}, and appear closer to those experiments than other methods of estimating it \cite{creechTearFilmThickness1998}.
\begin{figure}
    \centering
    \includegraphics[width=\textwidth]{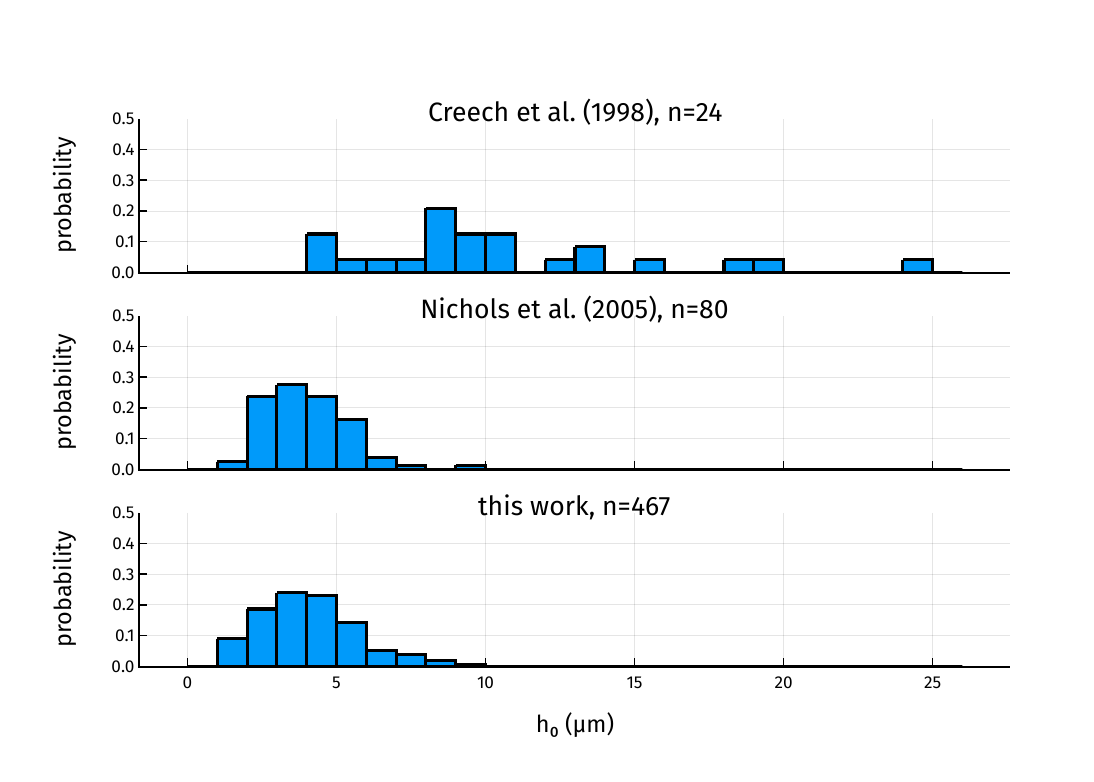}
    \caption{The probability distributions of initial TF thickness estimates from three sources:  Creech et al. \cite{creechTearFilmThickness1998} ($n=24)$,  Nichols et al. \cite{nicholsThinningRatePrecorneal2005} ($n=80$) using interferometry, and this work ($n=467$).
    }
    \label{fig:h0_hists}
\end{figure}
\figref{fig:h0_hists} shows histograms of the probability of the thickness from four sources, including the results from this work.  One can see that our pre-corneal tear film (PCTF) thickness estimates from fluorescein data agree well with the 80 interferometric measurements of Nichols et al. \cite{nicholsThinningRatePrecorneal2005}; that study discussed the possible sources of discrepancy with the relatively broad distribution of 20 no-lens estimates from Creech et al. \cite{creechTearFilmThickness1998} based on coating flow theory.  The 20 manual PCTF thickness estimates of Luke et al. \cite{lukeParameterEstimationMixedMechanism2021b} (not shown) form a narrow distribution that is easily within the experimental range \cite{nicholsThinningRatePrecorneal2005}.

\begin{figure}
    \centering
    \includegraphics[width=\textwidth]{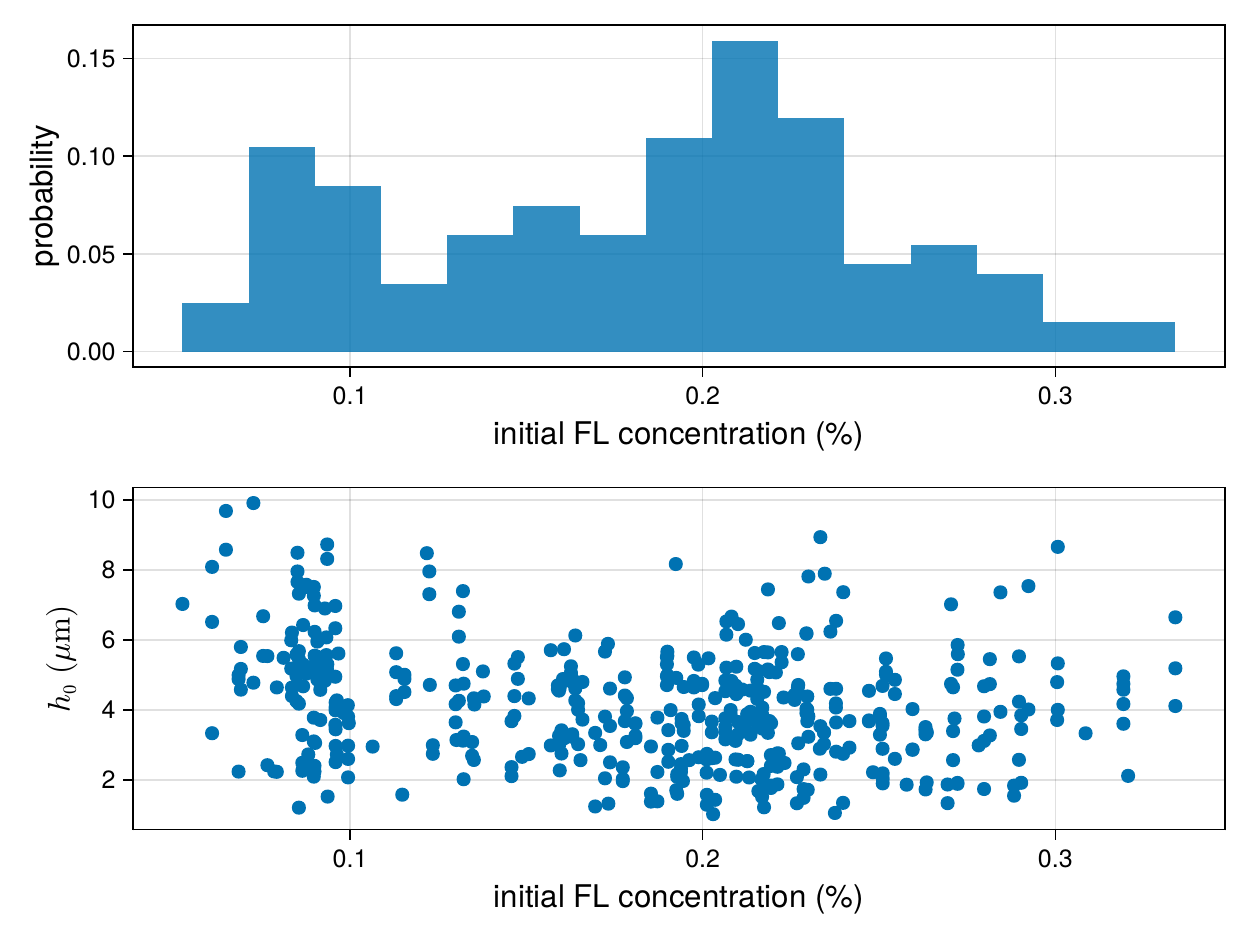}
    \caption{Top: Distribution by probability of estimated initial Fl concentrations $f_0$ over all trials. Bottom: Initial thickness $h_0$ vs $f_0$ over all fitted locations. The value of $f_0$ is the same for an entire trial, but $h_0$ can vary between locations within a trial.}
    \label{fig:f0_values}
\end{figure}

The initial thickness estimates require estimates of the initial fluorescein concentrations.
These are computed using the approach of Wu et al. \cite{wuIncreasingOSStimulation2015}.
The estimates for all subjects are shown in \figref{fig:f0_values}.
The upper part of the figure shows two peaks in the histogram.  Trials most often begin close to the critical concentration of 0.2 \%, which is the location of the right peak.  The left peak, around 0.1 \%, is due to the protocol for the experiments.  FL is not instilled for every trial so that one or two trials could occur before additional FL is instilled; tear turnover would reduce the FL concentration \cite{webberFLImaging1986}.  The lower part of the figure shows a scatter plot of the initial thickness and the initial FL concentration. The two do not appear to be correlated; the initial thicknesses seem uniformly spread across its range of values for all values of $f_0$.  As mentioned above, the distribution of thicknesses estimated from the $f_0$ in \figref{fig:h0_hists} agree quite well with measured distributions of thickness \cite{nicholsThinningRatePrecorneal2005}.

\subsection{Mechanism for all subjects}
\begin{figure}
    \centering
    \includegraphics[width=\textwidth]{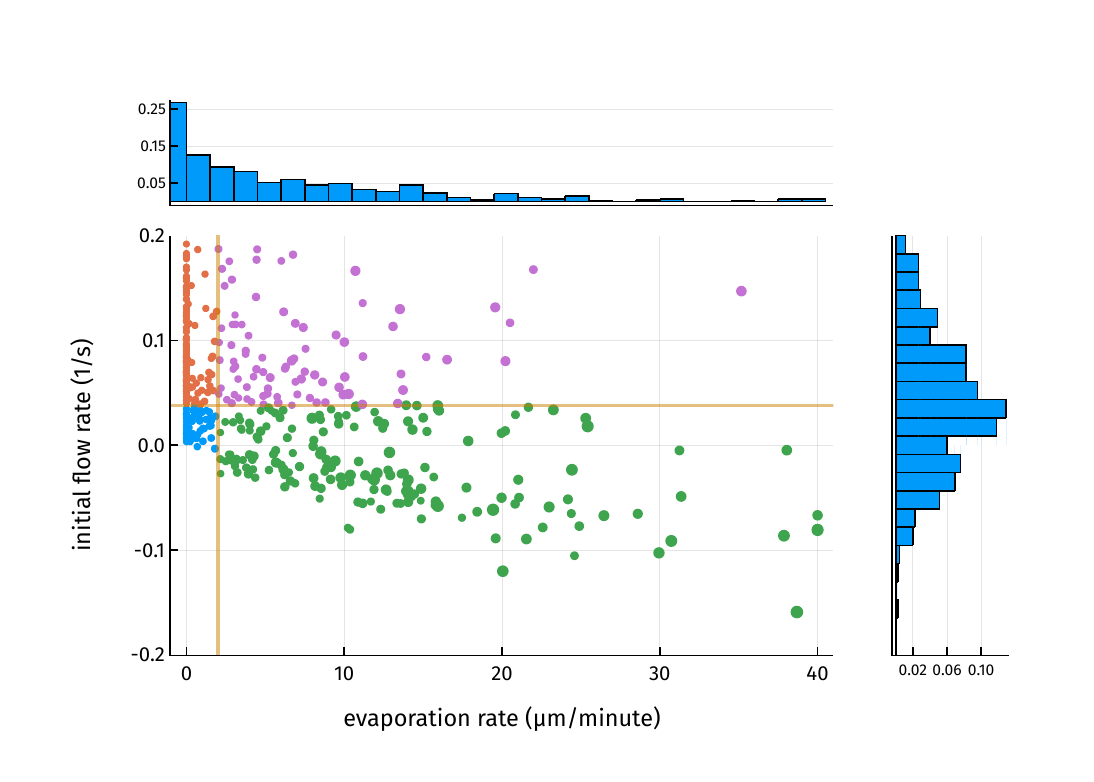}
    \caption{Scatter plot of evaporation and flow rates found from model-D fits to the data. The area of each dot is proportional to the final osmolarity predicted by the model. Marginal histograms show the distributions in probability of each parameter. The orange lines, drawn at \qty{2}{\mupermin} for evaporation and the median value \qty{0.038}{\per\second} for flow, are used to color each sample to indicate high/low rates of evaporation and flow.}
    \label{fig:evap_flow_osmolarity}
\end{figure}

\figref{fig:evap_flow_osmolarity} shows the results of all the fits as a scatter plot and marginal distributions of dimensional evaporation rate $v'$ and initial flow rate $b_1'$, with dot sizes indicating the final osmolarity value in the fitted model. The majority of the evaporation rates are at or below \qty{2}{\mupermin}, which agrees well with interferometric measurements of central cornea thinning rates \cite{nicholsThinningRatePrecorneal2005} and previous fitting work \cite{lukeParameterEstimationMixedMechanism2021b}. The specific choice of $2\mu$m/min was taken from the distribution of PCTF thinning rates in Nichols et al  \cite{nicholsThinningRatePrecorneal2005}, which showed a transition from a highly peaked set of low rates to a broad set of higher rates.

A key question for the models is the relative importance of evaporation rate and tangential flow.  We use the flow parameter $b_1'$, which is the initial strength of the tangential flow, as the indicator of the importance of flow. Large positive values indicate that divergent flow is important in tear thinning \cite{zhongMathematicalModelling2018,lukeParameterEstimationMixedMechanism2021a}, while a negative value indicates a convergent flow consistent with evaporation being a primary mechanism in thinning and TBU \cite{PengEtal2014,braunTearFilm2018}. \figref{fig:evap_flow_osmolarity} also shows lines drawn at \qty{2}{\mupermin} for evaporation and the median value \qty{0.038}{\per\second} for flow, partitioning the plot into four quadrants. The upper-left quadrant features TBU with low evaporation and strong divergent flow, suggesting that Marangoni-driven thinning dominates; a significant fraction of these cases seem to have little evaporation involved. The lower-right quadrant contains high-evaporation cases featuring little flow or an initially convergent flow whose strength generally increases with the evaporation rate. This scenario is consistent with inward tangential flow that tries to mitigate rapid evaporative loss \cite{PengEtal2014,braunTearFilm2018,lukeParameterEstimation2020}. The upper-right quadrant could be interpreted as mixed-mechanism, where both evaporation and outward tangential flow cooperate to thin the TF.  The lower-left quadrant represents cases that may not have enough thinning of either type to be definitively called TBU; we labeled these cases ``good tear film" or gtf.

\begin{figure}
    \centering
    \includegraphics[width=\textwidth]{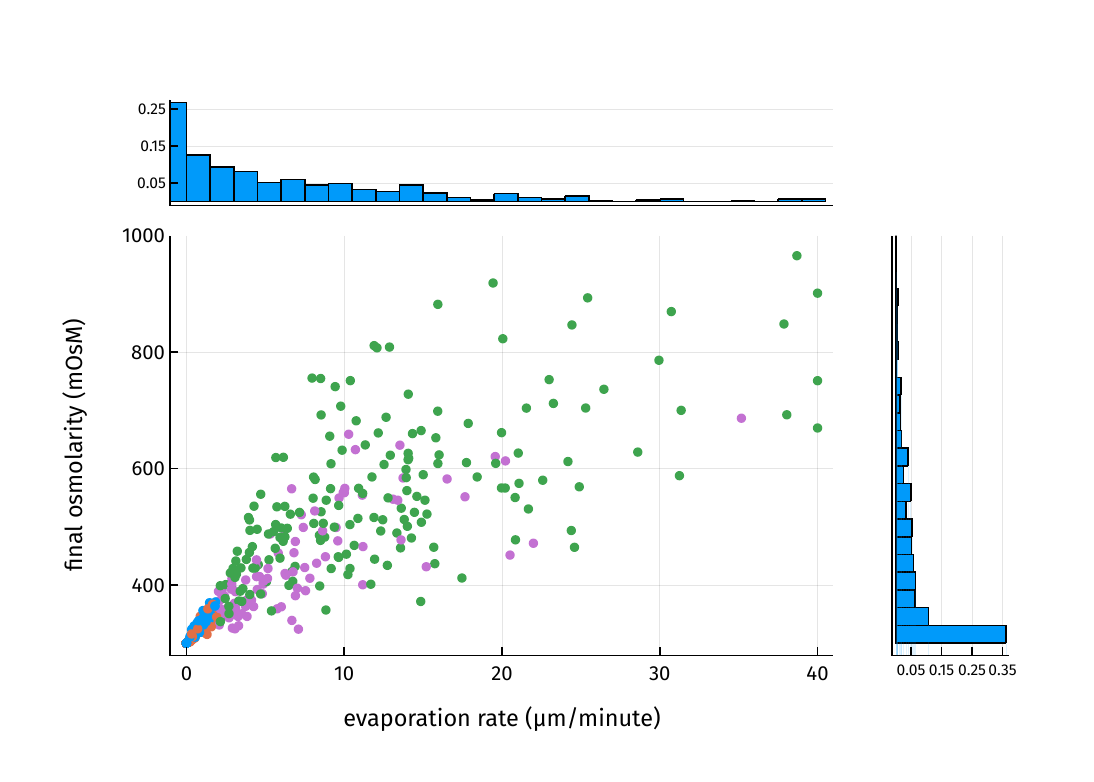}
    \caption{Scatter plot of evaporation rates $v$ and final osmolarities $c_e$ found from model D fits to the data for all subjects. Marginal histograms show the distributions in probability of the individual quantities. The coloring of the dots is the same as in \figref{fig:evap_flow_osmolarity}, indicating low/high values for evaporation rate and initial flow.}
    \label{fig:evap_osmolarity}
\end{figure}

It is clear from \figref{fig:evap_flow_osmolarity} that the osmolarity increases with the evaporation rate, and the relationship is plotted explicitly in \figref{fig:evap_osmolarity}. For reference, we note that normal tear film osmolarity measured from the inferior meniscus is somewhat variable \cite{TomlinsonTearFilm2006,JacobiTearFilm2011,gilbard1978}, reportedly averaging 301 mOsM in normal subjects, with diagnostic cutoffs for DED ranging from 305-318 mOsM \cite{TomlinsonTearFilm2006,JacobiTearFilm2011,LempTearOsmolarity2011,VersuraPerformanceTear2010}. However, a previous study suggests that the levels of tear film hyperosmolarity over the cornea could be as high as 800-900 mOsM \cite{liuLinkInstability2009}, much higher than the levels measured from the inferior meniscus \cite{JacobiTearFilm2011,LempTearOsmolarity2011}. In this study, the increase appears to be roughly linear for $v'$ below \qty{10}{\mupermin}, but then the osmolarity tends to level off and does not exceed 950 mOsM for this set of results. This trend agrees with previous fitting results on fewer TBU instances \cite{lukeParameterEstimationMixedMechanism2021a,lukeParameterEstimationMixedMechanism2021b} and models with strong outward flow \cite{zhongFLimaging2019}.
Previous theories of TF thinning and TBU that were not fit to experimental data could give higher final values of the osmolarity  \cite{PengEtal2014,BraunKing-Smith2015,braunTearFilm2018}.  The distribution of final osmolarity values shows that for these healthy subjects, the osmolarity remains below sensory threshold levels (450 mOsM\cite{liuLinkInstability2009}) in the majority of cases.

\begin{figure}
    \centering
    \includegraphics[width=\textwidth]{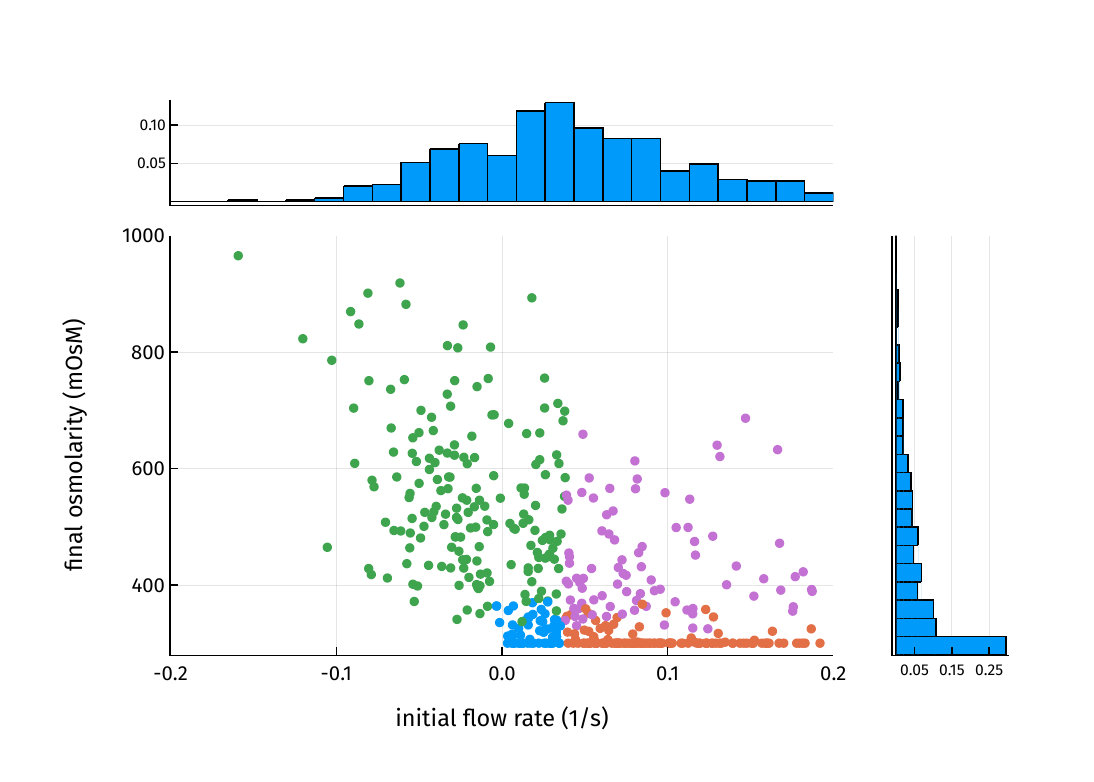}
    \caption{Scatter plot of initial flow rate $b_1$ and final relative osmolarity $c_e$ found from model D fits to the data for all subjects. Marginal histograms show the distributions in probability of the individual quantities. The coloring of the dots is the same as in \figref{fig:evap_flow_osmolarity}, indicating low/high values for evaporation rate and initial flow.}
    \label{fig:osmo_vs_b1}
\end{figure}

\figref{fig:osmo_vs_b1} shows a scatter plot of the final osmolarity vs the initial strain rate $b_1'$ for all subjects.  The final osmolarity is negatively correlated with flow:  many more hyperosmolar endpoints appear for low flow, and relatively few for stronger flow.  Fewer hyperosmolar endpoints at high flow may be expected, but the wide range of osmolarity that may occur at moderate or low flow is again apparent in these results.

The negative correlation observed for $v'$ and $b_1'$ help explain the results in the osmolarity.  There are relatively few cases where flow is important for larger $v'$, and so osmolarity is expected to become large in more of those cases.  What may be more surprising is that there is quite a range of flow strength for $2 \le v' \le 20\mu$m/min, and this causes a relatively wide range of values in the osmolarity for that range of $v'$.  There is an overall trend, but the flow and evaporation can cooperate to give high osmolarity in relatively short times, particularly for $v'\le 10 \mu$m/min.  This was seen in previous models to some degree \cite{lukeParameterEstimationMixedMechanism2021a,lukeParameterEstimationMixedMechanism2021b}, but with the current results this trend is more dramatic.  It is clear from these results that one cannot reliably estimate the final osmolarity from the TBU time alone; one needs knowledge of the local evaporation and flow conditions to get that estimate.

\begin{figure}
    \centering
    \includegraphics[width=\textwidth]{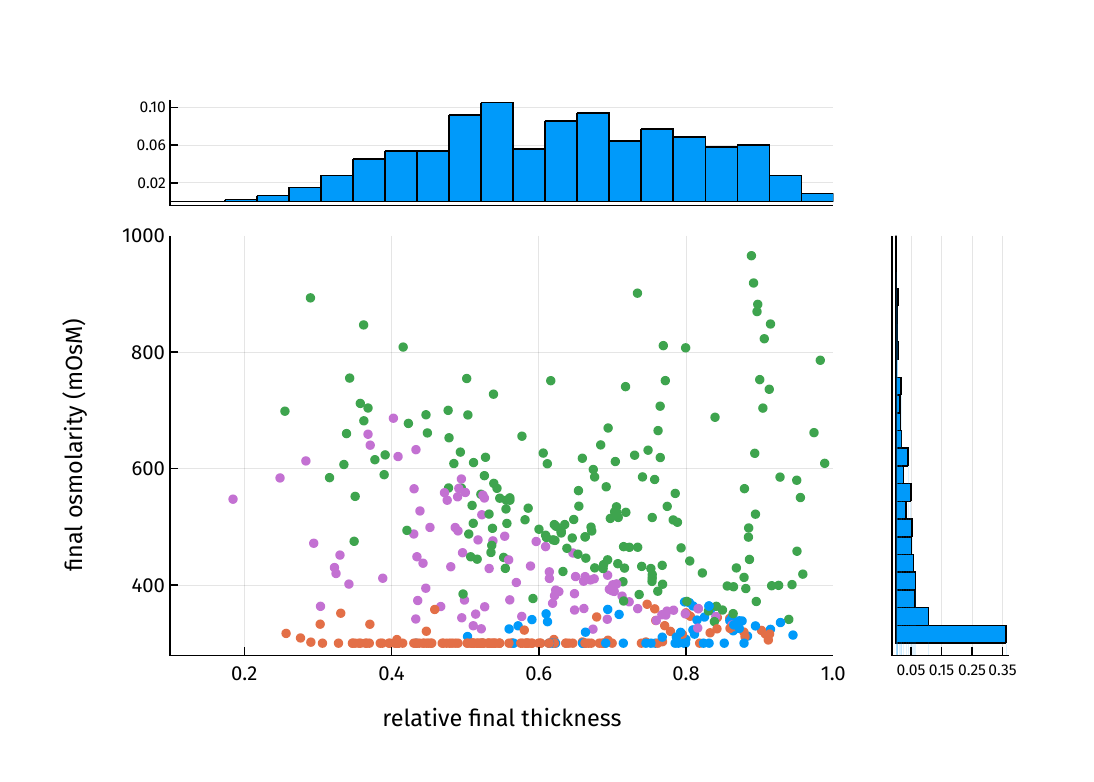}
    \caption{Scatter plot of relative final thicknesses $h_e/h_0$ and final osmolarities $c_e$ found from model D fits to the data for all subjects. Marginal histograms show the distributions in probability of the individual quantities. The coloring of the dots is the same as in \figref{fig:evap_flow_osmolarity}, indicating low/high values for evaporation rate and initial flow.}
    \label{fig:thickness_osmolarity}
\end{figure}

\figref{fig:thickness_osmolarity} shows a scatter plot and histograms for the final osmolarity $c_e$ and relative final thickness $h_e/h_0$. The high-flow cases (red and purple dots) are correlated with lower final osmolarity, but there is no clear association with the final thickness. We also note that the fit interval here may not extend to full thickness TBUT in many cases.  This is because the parameters can be strongly affected if the fit interval is too long or a late plateau of low intensity is included, and for that reason, the fit intervals could not be left too long.

\subsection{Evaporation and thinning comparisons}

\begin{figure}
    \centering
    \includegraphics[width=\textwidth]{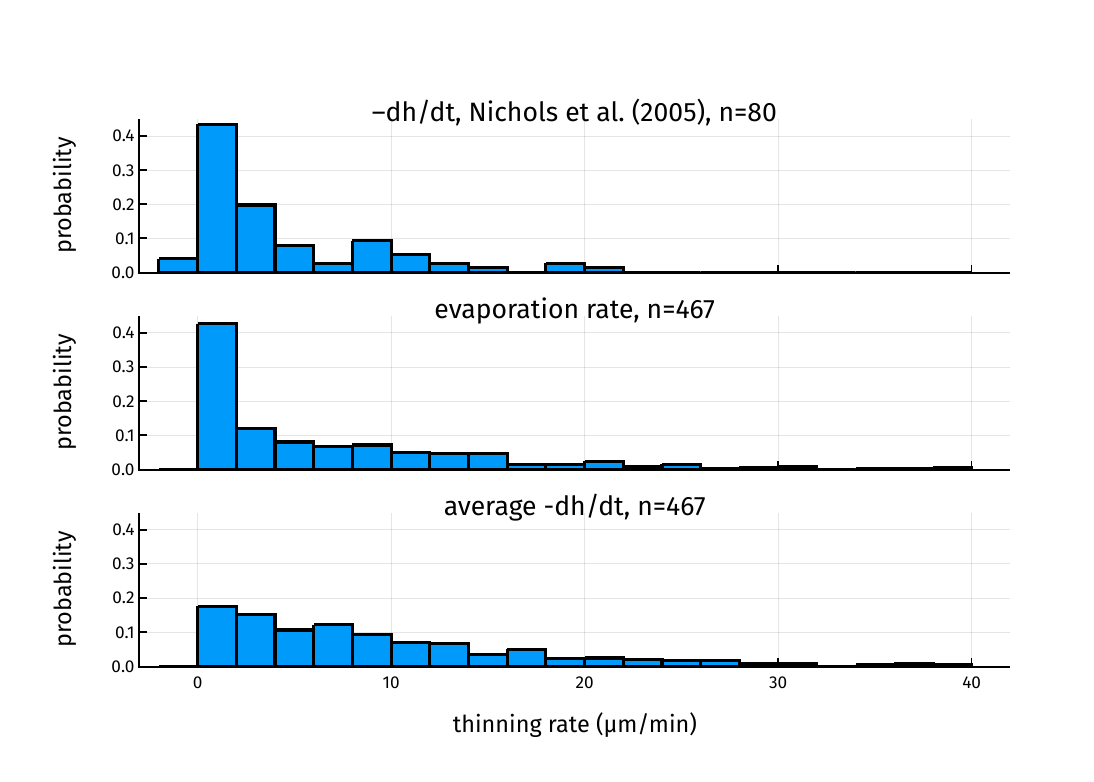}
    \caption{Comparison of thinning rate and evaporation  distributions from three sources
    experimental thinning rate $-dh/dt$ (with thickening in two cases) from Nichols et al. \cite{nicholsThinningRatePrecorneal2005}, as  measured in a central cornea spot of \qty{0.2}{\milli\meter} diameter; evaporation rate $v'$ from this work; and average $-dh/dt$ from this work.  See text for details.}
    \label{fig:evap_rates_compar}
\end{figure}

\figref{fig:evap_rates_compar} shows a comparison of thinning-rate and evaporation-rate results from the current work with experimental results~\cite{nicholsThinningRatePrecorneal2005}. The experiment used narrow-band interferometry to measure thickness rates centrally in a 0.2 mm diameter spot; thinning rates were computed from the slope of a best fit line that began 2 s after a blink.  The distribution of the measured thinning rates is within the values we found by fitting models to FL intensity decrease.  The evaporation rates we report do yield larger values than the thinning rates in experiment.  However, we note that the experiment can only detect intensity change with time, and it cannot separate the evaporative and flow-related contributions to TF thinning. Evaporative TBU can exhibit convergent tangential flow \cite{PengEtal2014,braunTearFilm2018}, which can significantly slow thinning.  Thus, it may be expected that evaporation rate could (and perhaps should) exceed the thinning rate in such cases.

\subsection{Results by subject}

\begin{figure}
    \centering
    \includegraphics[width=\textwidth]{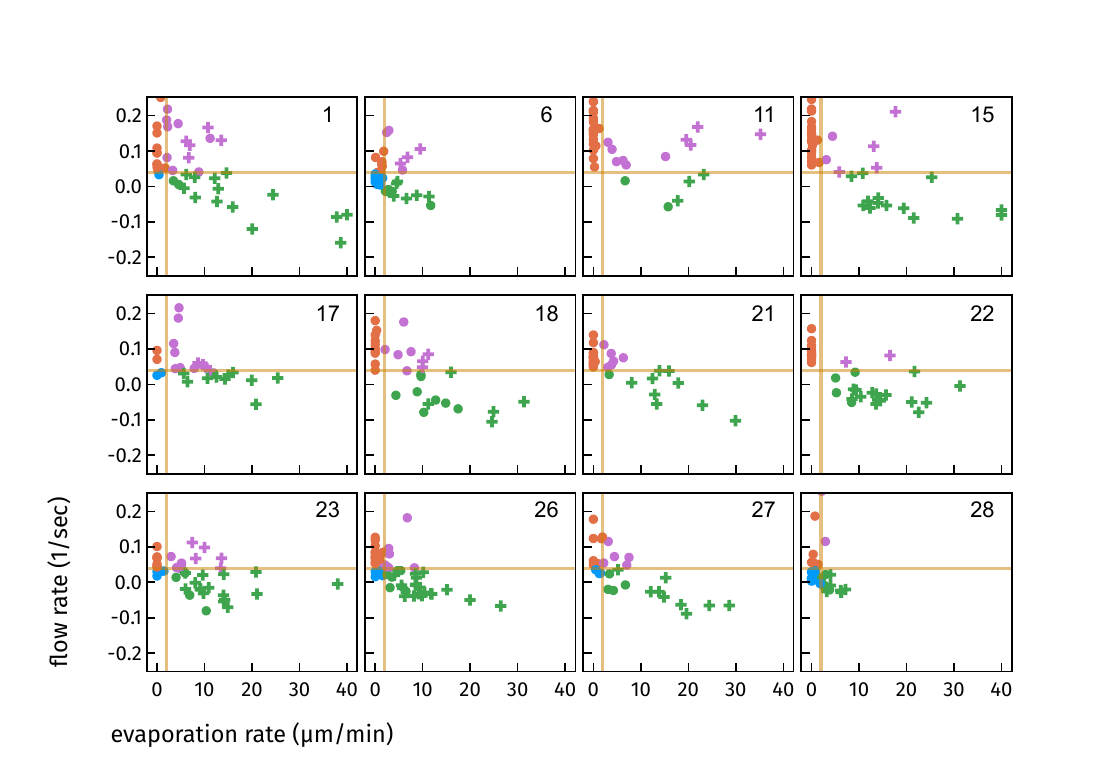}
    \parbox{6.25in}{\caption{Scatter plots of evaporation $v'$ (abscissa) and strain rates $b_1'$ (ordinate) separated by experimental subject. Three subjects who had fewer than 20 fits each have been omitted. Cases marked by a cross had a final osmolarity value above the 450 mOsM threshold of discomfort. The orange lines and symbol colors have the same meaning as in \figref{fig:evap_flow_osmolarity}.}
    \label{fig:scatter_by_subject} }
\end{figure}

\figref{fig:scatter_by_subject}, like \figref{fig:evap_flow_osmolarity}, shows the fitting results for evaporation rate $v'$ and initial strain (flow) rate $b_1'$, but plotted separately for each experimental subject. The plots also use a cross symbol to indicate the cases in which the final osmolarity exceeded the 450 mOsM discomfort threshold \cite{liuLinkInstability2009}. The orthogonal lines in each subplot show the boundaries that we chose between the different mechanisms.  High evaporation rate cases are to the right of the vertical line at $v'=2\mu$m/min; high flow cases are above the horizontal line at the median value of $b_1'$, which is 0.038 s$^{-1}$ (computed over all trials and instances).  The coloring scheme for the symbols is the same as in \figref{fig:evap_flow_osmolarity}.

No subject displayed exclusively high-evaporation TBU, while a few were characterized by low evaporation rates. Most of the low-evaporation, low-flow cases that may indicate good TF occurred in just three subjects.
All subjects displayed some cases of both positive (divergent) and negative (convergent) initial transverse flow, at rates dispersed rather widely in most cases. It was fairly common to exceed the discomfort threshold due to high evaporation marked by convergent flow, while it was less common to exceed the threshold with high divergent flow.

Table~\ref{tab:all_subjects_mech} shows results for the healthy subjects we studied.  The table is in descending order of number of instances fit for each subject.  The four possible mechanisms are that (i) evaporation drives TBU; (ii) flow drives TBU; (iii) a mix of evaporation and flow drives TBU; and a ``good TF" (gtf) where neither evaporation nor flow is very strong.
The results show that on a population level, evaporation is the most common driver of TBU.  The second most common instance is flow, which has both small evaporation and relatively large flow rates.
The relative position with respect to $v'=$ 2\micron/min and $b_1'= 0.038$ s$^{-1}$ divided instances into these four classes; if greater than these threshold values, the effect was important, and vice versa if less.

\begin{table}
\begin{centering}
\begin{tabular}{|c|c|c|c|c|c|}
  \hline
  subject & evap & flow & mixed & gtf & total \\ \hline
  26 & 25 & 14 & 7 & 5 & 51 \\
  15 & 14 & 24 & 6 & 0 & 44 \\
  1  & 16 & 12 & 14 & 1 & 43 \\
  23 & 20 & 7 & 10 & 5 & 42 \\
   6 & 12 & 8 & 7 & 13 & 40 \\
  28 & 13 & 6 & 2 & 17 & 38 \\
  22 & 20 & 13 & 2 & 0 & 35 \\
  11 & 5 & 18 & 10 & 0 & 33 \\
  21 & 10 & 16 & 6 & 0 & 32 \\
  18 & 13 & 10 & 8 & 0 & 31 \\
  27 & 13 & 7 & 6 & 3 & 29 \\
  17 & 11 & 2 & 11 & 2 & 26 \\  \hline
  2 & 4 & 3 & 0 & 10 & 17 \\
  10 & 0 & 0 & 2 & 1 & 3 \\
  14 & 0 & 1 & 2 & 0 & 3 \\  \hline
  all & 176 & 141 & 93 & 57 & 467 \\
  \hline
\end{tabular}
\caption{For the 15 subjects that were fit, there were four possible mechanisms: evaporative (evap), flow, mixed (evaporation and flow) and good TF (gtf).  The distribution of mechanisms for all instances fit is given for each subject in descending order of number of instances fit.  We excluded the last three subjects from any subject-specific analysis because of the few instances of TBU fit.}  \label{tab:all_subjects_mech}
\end{centering}
\end{table}

We found that perturbing these threshold values around these points did not affect the relative distribution of the mechanisms very strongly.  However, there is still some dependence on the choice of the boundaries.  For example, if we used the overall median values for both $v'$ and $b_1'$, then the number of evaporation driven and flow driven instances are equal, and the number of mixed and gtf cases are equal but at about half the number of the other categories.  While using the median clearly has statistical rationale, we found experimental motivation for our choice of a different $v'$.  In contrast, there is no experimental guidance for $b_1'$ since there are no direct measurements of flow in thinning, so we used this statistically-based choice for this parameter.

From Table~\ref{tab:all_subjects_mech} we also see that the distribution of mechanisms may be different from subject to subject.  For example, subject 26 has a preponderance of evaporative TBU cases, while subject 15 has a preponderance of flow TBU mechanisms.  Furthermore, subjects 6 and 28 have a large fraction of gtf cases.  This suggests that, in some cases, it is possible to distinguish subjects based on their TBU mechanism distribution.  Though the specifics numbers may change for different mechanism selection criteria, the distribution of values will still typically vary from subject to subject.

\begin{figure}
     \centering
     \includegraphics[width=\textwidth]{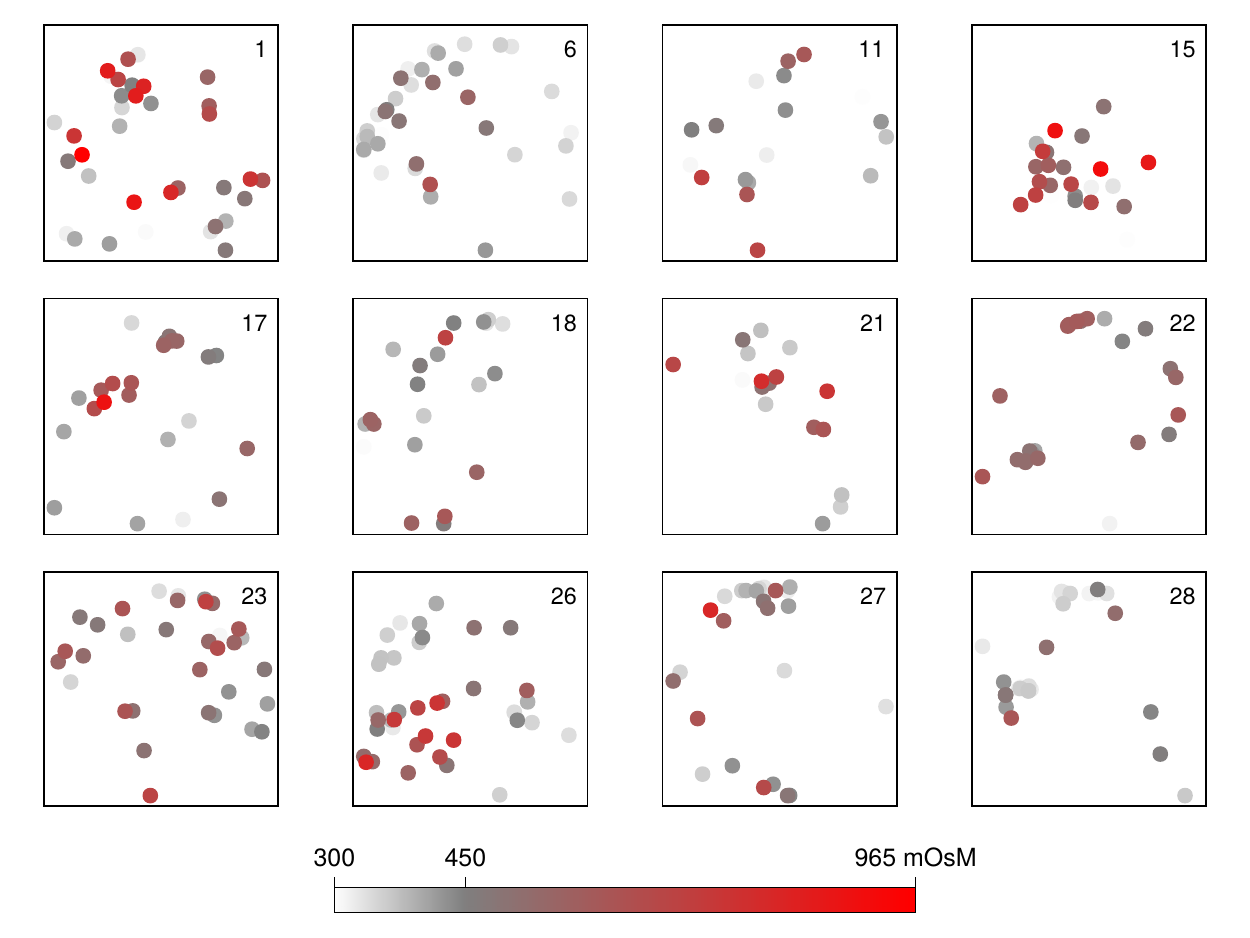}
     \caption{Scatter plots showing the locations of every fitted TBU for each of the top 12 subjects. Each spot is located relative to a box fitted approximately around the cornea throughout the trial; the data and squares are shifted and stretched slightly for alignment in this figure, and the position of the cornea within the box is not constant during a trial or between different trials. The color shows the model's final value for osmolarity, with gray for values below the discomfort threshold of 450 mOsM \cite{liuLinkInstability2009} and increasing saturation of red for values above the threshold.}
     \label{fig:osmolaritymap_by_subject}
\end{figure}

\figref{fig:osmolaritymap_by_subject} uses scatter plots to show the locations of every fitted TBU instance for each of the 12 experimental subjects with at least 20 total fits. The location of each dot is relative to a box approximately enclosing the average detected cornea throughout its experimental trial; however, the cornea does not have fixed position within the box during a trial or between different trials. The color of each dot is used to indicate the fitted model's final value of the osmolarity, with cases above the discomfort threshold colored red. Each subject exhibited some potentially painful TBU instances, which is consistent with the instructions given to the subjects in the experimental trials.  Because this data came from sustained tear exposure (STARE) trials, this distribution represents the stimulus to the cornea in life outside the clinic.

\section{Discussion}
\label{sec:discuss}

In this paper, we have employed established methods of fitting mathematical models to FL intensity data \cite{lukeParameterEstimationMixedMechanism2021b} and applied them to automatically identified thinning or TBU regions in healthy subjects.  We modified the approach of Su et al. \cite{suTearFilm2018} to identify multiple TBU regions of interest in each trial.  With those automatically identified thinning regions, we had 12 subjects with a significant number of TBU instances over 2 visits and 10 trials (though no subject yielded fits from all trials).  The mechanisms for each thinning and TBU is determined from the optimal coefficients from fitting the FL intensity data.   We find that each subject has a range of mechanisms associated with
their sample of TBU instances, and those distributions of mechanism may be sufficient to distinguish between subjects.  By pooling the results from all of the subjects, we have a relatively large sample of 467 instances of thinning or TBU from 15 subjects.

A primary result is that the final osmolarity varies widely across the set of thinning or TBU instances that we studied.  This appears to happen for most subjects individually as well as for the pooled instances.  For very rapid thinning with little evaporation, or instances that turn out to thin very little, the final osmolarity stays low.  For instances that are driven by evaporation or mixed mechanism, the final osmolarity rises to as much as  750 mOsM except for uncommon values reaching over 900 mOsM.  Because these results vary by subject, the results have the potential to distinguish between subjects.

The results in this paper expand on prior efforts to fit FL intensity data in TF thinning and TBU.   PDE models showed the potential to fit the intensity TBU instances in some cases \cite{braunTearFilm2018,zhongMathematicalModelling2018}, and to determine mechanisms that drive thinning and TBU \cite{lukeParameterEstimation2020,lukeParameterEstimationMixedMechanism2021a}.
Simplifying to local ODE models of TBU allowed for faster computation and freedom to select more instances of TBU; 20 instances were fit by ODE models in Luke et al \cite{lukeParameterEstimationMixedMechanism2021b}.  The data in this paper confirmed previously observed trends \cite{lukeParameterEstimationMixedMechanism2021b} with more than 20 times the instances fit.  We are unaware of other work that fits the dynamics of TBU with mathematical models or estimates the parameters within TBU at this scale.

Experimental measurements of initial tear thicknesses include a range of approximately 2 to 10 $\mu$m using interferometry \cite{King-SmithFink2004,nicholsThinningRatePrecorneal2005} and similar values from optical coherence tomography (OCT) \cite{wangTFThicknessOCT2003}.  The distribution of our initial thickness estimates matches those of interferometry measurements very well. The values found via estimation in Creech et al. \cite{creechTearFilmThickness1998} give a wider range of thicknesses and a substantial fraction of thicker film estimates.
Dursch et al. \cite{DurschFLandThermal2017} fit a model to the thinning and the temperature of the TF to estimate the evaporation rate of the TF.   To our knowledge, they did not determine the osmolarity of the TF in TBU, but they did use imaging data from both FL intensity and thermal imaging to determine TF parameters of interest.

We now turn to the strengths and weaknesses of our approach.  The method identifies multiple instances of thinning and/or TBU from each trial.  The system uses a trained CNN to find regions of interest from which minimum values of intensity in the ROI are extracted throughout the trial.  The training of the CNN used labeled images with a fixed threshold of intensity for TBU across all trials; the overall intensity of the trials varied, however, so it would likely expand the number of trials that could be analyzed to make that threshold for TBU trial dependent.  The ROIs are found near the end of the trial, and this approach assumes that the thin regions are not moved around by flow.  However, it is possible that thinning begins elsewhere and flow moves the thinning spot into the ROI during the first seconds of the trial \cite{King-SmithIOVS13a}; this type of dynamic is beyond what our model can analyze at this time.  Extracting the minimum intensity additional Gaussian blurring in the ROI gives acceptable FL intensity data for fitting; however, the data is noisy even after the Gaussian filtering and smoothing.  There are other possible choices of method for extracting the intensity data, but our approach did not seem to be too sensitive to what we attempted.
Some instances have rather little happening, but are relatively dark compared to their surroundings.  One could ask whether any ``good TFs" should have been selected for fitting.  It is unclear at the time of writing whether this should be the case or not; a rather long trial looks like nothing is happening but eventually TBU may occur, yet would still be rated as a good TF.  Not all trials or subjects can be analyzed by the automatic system for TBU.  This may be caused by lack of an inferior meniscus for estimating FL concentration, failure of the initial thickness estimate, poor recording of intensity data due to subject movement or to image focus, or possibly other reasons. The relatively low number of subjects could also be considered a limitation of the study.

Fluorescein is used for imaging which may impact TF dynamics \cite{mengherEffectFLonTBU1985,ChoDouthwaite95,mooiMinimizeFL2017}.   The initial FL concentration estimates show some systematic variation, possibly because two to three trials were performed between each FL installation. It is not clear whether a different installation protocol may improve our method.  The variation in initial FL concentration did not pose any difficulty for estimating the initial thickness in the cases that were fit by the models, and may reduce variation between visits that may affect other tests such as TBUT determination \cite{coxNIKBUTvsTBUT2015}. The use of STARE trials does not represent healthy blinking but it does ensure that thinning and TBU occurs.

Despite these limitations, the method has produced repeatable data for hundreds of instances of thinning and/or TBU.  The data reveals trends in the conditions experienced by a cohort of healthy TF subjects.
According to the model, within TBU the final osmolarity is highly variable due to the differing mechanisms driving TBU;  this is lower than the upper limit suggested by some previous models without fitting \cite{PengEtal2014,BraunKing-Smith2015,braunTearFilm2018}, higher than flow driven models initially suggested \cite{zhongMathematicalModelling2018} and agrees with previous models that fit FL intensity \cite{lukeParameterEstimationMixedMechanism2021b}.
The final osmolarity may be high within TBU but appears to stay below 950 mOsM for this set of subjects, in agreement with the result of Liu et al \cite{liuLinkInstability2009}.   The evaporation and thinning rates appear to agree well with published data \cite{nicholsThinningRatePrecorneal2005}, and
the relationship between evaporation rate and final osmolarity is revealed to be generally increasing with evaporation rate but is complicated by the dependence on flow.
The model determines optimal flow and evaporation values and the direction of flow (from the sign of $b_1'$) is a major part of determining the mechanism of an instance.

The DEWS II Diagnostic Methodology report \cite{wolffsohnTFOSDEWS2017} recommended using non-invasive TBUT to help diagnose DED rather than fluorescein due the variations induced by the latter \cite{mengherEffectFLonTBU1985,choReliabilityTear1992,choTearBreakup1998}.
The utility  of each type of method is still an active area of research \cite{mooiMinimizeFL2017,paughDEDandTBUTmethods2019,speakmanDiagnosticUtilityNITBUT2022}.
Many of these approaches average the results of two or three measurements, and may eliminate outlying values, as suggested by Cho et al.\ \cite{choReliabilityTear1992} and many others. Recent efforts have tried to automate TBUT determination  \cite{suTearFilm2018} and DED diagnosis \cite{remeseiroCASDESComputerAided2016}. The study of Segev et al.\ \cite{SegevDynamicAssessment2020} found breakup times based on mean values of aqueous layer thickness from two 40s trials separated by 45 min on average. We note that our approach is not aimed at using TBUT as a method for diagnosis; we are refining the use of FL imaging to yield the mechanism driving thinning and TBU for \emph{many} instances in each healthy subject.  We are attempting to find the distribution of what can occur within healthy subjects, and there appears to be significant variability within each subject and even within a single trial. We are unaware of prior studies that investigated within-subject TBU variability.  This basic science data may have clinical application in classifying subjects based on their thinning and/or TBU characteristics.

Some studies have noted and tried to exploit the distributions of tear film parameters to distinguish between subjects.  An example is Bai et al \cite{baiTFLLwithOptical} where optical microscopy is used to measure the LL thickness for healthy subjects and several conditions related to meibomian gland dysfunction.  The distribution of LL thickness over a small area is analyzed for each subject and differences between conditions can be seen from these distributions.

In this study, we identify parameters and mechanisms for multiple instances of thinning in each subject.  Those instances present varying amounts of chemical, thermal and mechanical stimuli to the ocular surface.  The mechanism by which those stimuli are sensed or received, and the role of that perception in DED, is a matter of ongoing research \cite{dewsIIPain}.   Various neural receptors are thought to play important roles in sensing these different stimuli:  chemical \cite{dewsIIPain}, thermal \cite{parraHyperosmAndCold2014,hirataHyperosmAndCold2014,situTFInstabilityResponses2019}, and mechanical \cite{awisi-gyauChangesCorneal2019,belmontePainDrynessItch2019}. While our work here cannot directly address such questions, we believe that quantifying the stimulus at the ocular surface can only help to clarify such processes.

In order to compute the fits to the extracted data, reasonable ROIs for extraction must be found in each trial. For the healthy subjects that we used, less than half of the trials yielded ROIs for analysis.  Improving the robustness of the ROI detection would be an efficient way to generate more data to characterize thinning and TBU instances and the subjects in which they occur.  Once ROIs are determined, there may be other options than what we employed for extracting the thinning data.  The FL images were somewhat noisy, and despite filtering to minimize it, extracting local data may be affected by that noise.

The estimates of initial FL concentration, and subsequently the initial thickness, required a special procedure with low illumination intensity and a good inferior meniscus.  This may limit the trials and subjects that may be analyzed.  Other possible ways to estimate these initial quantities may improve robustness of the method.  The method appears to work very well when estimates can be obtained, based on comparison with interferometric \emph{in vivo} results.

\section{Conclusion and future perspectives}

An important next step would be to apply the method to a sample of DED subjects to compare with the data from healthy subjects. Combining our method with data from simultaneous thermal imaging\cite{DurschFLandThermal2017}, interferometry\cite{King-SmithIOVS13a}, or sensory feedback and/or sensory response \cite{liuLinkInstability2009,awisi-gyauChangesCorneal2019,liPainLinkIBI2018} could yield new insights.

\appendix
\section{Model Derivation}

Consider a rectangular control volume of $-L'/2 \le x' \le L'/2$ and $0\le y' \le h'(t')$; this rectangle could be centered on \figref{fig:model_sketch}.  The equations result from conserving solvent (the aqueous layer's water) and solutes (osmolarity, $c'$, and fluorescein concentration, $f'$, both in M) per unit width of the film. The water conservation is given by
\begin{equation}
\rho L'\frac{dh'}{dt'} = -J_e'L'+\rho P_o V_w(c'-c_0) L' -2 \rho h'  u'(L'/2,t')
\end{equation}
where the (constant) evaporation rate is given by $J_e' = \rho v'$.
The term on the left is the rate of change of the mass of water in the control volume.  The first term on the right is the water lost due to evaporation; the second term is supply of water due to osmosis.  The remaining term is the total amount of water flowing out of the ends from a depth-independent velocity field $u'(x',t')=g'(t')x'$; this velocity along the film is evaluated at $x'=\pm L'/2$.  The time dependence is given by $g'(t')=b_1'e^{-b_2't'}$.

Conservation of solutes is given by
\begin{equation}
   L'\frac{d(h's')}{dt'} = -[ 2 u'(L'/2,t') ] h' s',
\end{equation}
where $s'= c'$ or $f'$.  The term on the left is the rate of change of solute in
inside the control volume, and the term on the right is the total amount of solute leaving the sides of the control volume at $x'=\pm L'/2$.

Substituting $u'(L'/2,t')=g'(t')L'/2$ into the equations and rearranging gives, for water,
\begin{equation}
\frac{dh'}{dt'} = -v' + P_o V_w(c'-c_0)  - h'g'(t'),
\end{equation}
and for solutes,
\begin{equation}
   \frac{d(h's')}{dt'} = - g'(t') h' s'.
\end{equation}
Substituting for $g'(t')$ and converting to non-dimensional variables via (\ref{eq:hc_scales}) and $f'=f_{cr}f$ results in the nondimensional equations given in
Section~\ref{sec:modelsIntro}.

\section{Physical parameter values}

\begin{table}[htb]
  \begin{center}
  \caption{Dimensional parameters. Molar extinction coefficient\cite{mota1991} has been multiplied by ln(10) to convert it to the Napierian form.}
   \label{tab:dim}
\def~{\hphantom{0}}
\begin{tabular}{cll}
Parameter & Description & Value \\
$\mu$      & Viscosity \cite{tiffany1991}                            & 1.3$\times10^{-3}$Pa$\cdot$s            \\
$\sigma_0$ & Surface tension \cite{nagyova1999}                      & 0.045N$\cdot$m$^{-1}$                    \\
$\rho$     & Density (water)                             & $10^3$kg$\cdot$m$^{-3}$                \\
$h_0$       & Characteristic thickness\cite{King-SmithFink2004}             & $1$ to $10\ \mu$m              \\
$t_s$   & Characteristic time scale           & [fit interval in s]                      \\
$P_o$  & Tissue permeability of cornea  \cite{BraunKing-Smith2015}                        &  12.0$\mu$m/s \\
$V_w$      & Molar volume (water) &  $1.8\times 10^{-5}$m$^3\cdot$mol$^{-1}$       \\
$\epsilon_f$   & Naperian extinction coefficient\cite{mota1991}      & $1.75 \times 10^{7} $ m$^{-1}$M$^{-1}$   \\
$c'_{0}$   & Isotonic osmolarity \cite{gilbard1978}                 & $302$ mOsM$^3$                        \\
$f_{cr}$   & Critical fluorescein concentration \cite{nichols2012}   & $0.2\%$ (by mass)                       \\
$v$ & Experimental thinning rate\cite{King-SmithFink2004,WinterAnderson10} & $-3$ to $25$ $\mu$m/min 
  \end{tabular}
  \end{center}
\end{table}

The dimensional parameter values are given in Table~\ref{tab:dim}. Dimensionless parameters using typical values are given in Table~\ref{tab:nondim}.  Note that because $\phi$ varies by instance of thinning, and $P_c$ varies in the optimization to fit each instance of thinning, we only give typical values here.
\begin{table}
\begin{center}
\caption{Dimensionless parameters that arise from scaling the dimensional fluid mechanics problem. The values given are based upon the values of Table \ref{tab:dim}, $h_0 = 3 \  \mu$m, and $t_s = 3$ s.}
\label{tab:nondim}
\begin{tabular}{cccc}
 Parameter &  Description   &  Expression   &  Value  \\  
$P_c$   & Permeability of cornea    & $\displaystyle P_o V_w c_0/(h_0/t_s)$              & 0.0653  \\
$\phi$   & Napierian extinction coefficient   & $\displaystyle \epsilon_f f_{\text{cr}} h_0$                & 0.279  \\ 
\end{tabular}

\end{center}
\end{table}

\section{Computational details}

Here we give some details of the numerical procedures used to identify TBU instances and to process and extract the intensity data.

The CNN described in section~\ref{sec:tbu_detect} was trained in TensorFlow version 2. The loss function used was categorical crossentropy, the optimizer was ADAM, for metrics we used `accuracy' and the batch size was 32.

To apply the trained CNN to find instances of thinning and TBU, each image is cropped close to a circle automatically fit closely to the cornea \cite{driscollAutomaticDetectionCornea2021}. The CNN predictor is applied within a $192\times192$ window that is moved in overlapping fashion using a stride of 32 pixels. The window locations having a predicted TBU probability of 0.999 or higher are collected and clustered with an unsupervised hierarchical clustering algorithm (\texttt{hcluster} in scipy with 75 as the distance criterion).  If closer than the criterion, clustered tiles are merged.  The process is started from the beginning of a trial and continues until at least three distinct instances of apparent TBU were found.  After this process, a $192\times192$ box is centered on each cluster to serve as the TBU ROI.

Once the ROI is found, the portion of each image in the ROI box is downsampled to $96\times96$, converted to gray scale, and blurred using a $21\times21$ Gaussian filter with standard deviation $\sigma=5$.  Then the minimum pixel location for each frame in the extraction interval was found.  The location to extract the intensity data was by taking the median of the minimum locations over all images after the first three seconds of the trial.

Prior to analyzing the extracted data, three indicators of time series quality are computed.  First, the fraction of ROI minimum points lying within a $33\times33$ window around the extraction location serves as a check on how much variance was in the minimizing point over trial time. Second, if the pixel intensity at the extraction location when TBU is first predicted is greater than the value of 60 used for training the CNN, then that instance might not truly represent TBU. Third, if the median intensity over the first 10\% of the time series compared to the median intensity over the last 10\% did not indicate a decrease of at least 25\%, the time series might not show true thinning. Any ROI that raised an exception to the tests was checked manually for inclusion or exclusion.

Fitting the data to models was performed in Julia. Numerical solution of an ODE model with proposed parameter values is computed by the \texttt{DifferentialEquations} package with relative and absolute error tolerances of $10^{-10}$ and $10^{-11}$, respectively. A solution is immediately terminated with a large penalty if either $I(t)$ or $h(t)$ is found to be instantaneously increasing. The misfit of a proposed numerical solution is calculated as the trapezoidal 2-norm of the difference between numerical solution and data. The misfit is minimized by the \texttt{NLopt} package with box constraints and relative and absolute error tolerances of $10^{-5}$ and $10^{-7}$, respectively. The optimizer is initialized and run multiple times, with more complex model types including the values found by simpler models for the same data.  Optimization is performed by Nelder--Mead but was checked by Levenberg--Marquardt for consistency over all the time series.

Further details about the solution process can be obtained by inspecting the code in the repository at \url{https://github.com/tobydriscoll/fitting-ode-models-tear-film-breakup}.

\section*{Acknowledgements}

This work was supported by National Science Foundation grant DMS 1909846. The content is solely the responsibility of the authors and does not necessarily represent the official views of the funding source.

The authors thank P.\ Ewen King-Smith for helpful conversations and for sharing data for thickness and thinning rates.

\printbibliography

\end{document}